\documentclass{amsproc}
\usepackage{amssymb,amsmath,amsthm}
\usepackage{pdfpages,hyperref}

\newcommand\ovl[1]{\overline{#1}}
\newcommand{\gen}[1]{\langle #1 \rangle}
\newcommand{\eps}{\varepsilon}
\newcommand{\Z}{\mathbb{Z}}
\newcommand{\C}{\mathbb{C}}
\def\P{\mathbb{P}}
\newcommand{\izo}{\simeq}
\newcommand{\ra}{\rightarrow}
 \newtheorem{thm}{Theorem}[section]
 
 \newtheorem{lem}[thm]{Lemma}
 \newtheorem{prop}[thm]{Proposition}
 \theoremstyle{definition}
 \newtheorem{defn}[thm]{Definition}
 \theoremstyle{remark}
 \newtheorem{rem}[thm]{Remark}
 \numberwithin{equation}{subsection}
\begin{document}

\title[On Kummer 3-folds]
 {On Kummer 3-folds}

\author{ Maria Donten }

\address{Instytut Matematyki UW, Banacha 2, PL-02097 Warszawa}

\email{marysia@mimuw.edu.pl}





\begin{abstract}
We investigate a generalization of Kummer construction, as introduced in \cite{AW}. The aim of this work is to classify 3-dimensional Kummer varieties by computing their Poincar\'e polynomials.
\end{abstract}

\maketitle

\section*{Introduction}
A recent paper by Andreatta and Wi\'sniewski \cite{AW} provides a description of a generalization of Kummer construction, which is a method of producing a variety by resolving singularities in a quotient of a product of abelian varieties by a finite integral matrix group action. Some restrictive assumptions, both on the group action and the resolution, assure that a variety $X$ obtained as a result of the Kummer construction is projective, has $K_X$ linearly trivial and $H^1(X, \Z) = 0$.

In this paper the Kummer construction is applied to the product of three elliptic curves $A^3$. We look at the action of a finite group $G < SL(3,\Z)$ on $A^3$, which come from the natural action on $\Z^3$ by the identification $A^3 = \Z^3 \otimes_{\Z}A$. The quotient $A^3/G$ is singular. By resolving its singularities we obtain a Kummer 3-fold.

An important observation is that 3-dimensional Kummer varieties are Calabi--Yau.

The aim of this work is to compute Poincar\'e polynomials of such Kummer 3-folds, i.e. polynomials $P_X(t) = \sum_{i=0}^{2n} b_i(X)t^i$ of a formal variable $t$ with $b_i(X)$ being the $i$-th Betti number of a variety $X$. By $H^i(X, \C)$ or $H^i(X)$ we denote the complexified De Rham cohomology of $X$.

Section \ref{finite_subgroups_SL3Z} discusses the classification of finite subgroups of $SL(3, \Z)$. We describe the process of determining all groups to which we apply the construction. They are listed in \ref{groups_table}. In section \ref{construction_details} the details of Kummer construction in 3-dimensional case are explained. We concentrate on understanding the structure of singularities of the quotient of $A^3/G$ and their resolution. Thus we explain the method of computing Poincar\'e polynomials of Kummer $3$-folds. The remaining sections are devoted to the presentation of the results of our computations. All interesting points are clearly visible in the examples \ref{case_S4_2} and \ref{case_S4_3}, which are therefore discussed in details. In the remaining cases (in section 6.) only the result and most important data about the group action is given. The results of the computations are summarized in \ref{summary}, theorem \ref{last_theorem}.

The following results constitute the main part of my M.Sc. thesis completed at the University of Warsaw. I am greatly indebted to Jaros\l{}aw Wi\'sniewski, my thesis advisor, and I wish to thank him for all his help during the preparation of this paper.


\section{Finite subgroups of $SL(3, \Z)$}\label{finite_subgroups_SL3Z}

The first step is determining groups to which we apply the construction, i.e. finite subgroups of $SL(3, \Z)$. We say that $G, H < SL(n, \Z)$ are $\Z$-equivalent if they are conjugate in $GL(n, \Z)$.

\begin{rem}\label{conjugacy_classes}
The Kummer construction produces isomorphic varieties for groups which are $\Z$-equivalent.
\end{rem}

This means that we only need to find representatives of $\Z$-equivalence classes, called $\Z$-classes for short. We will restrict our attention to non-cyclic subgroups, because by lemma 3.3 in \cite{AW} cyclic groups fail to satisfy the assumptions of the construction.

We were not able to find the reference for the classification of finite subgroups of~$SL(3, \Z)$ up to $\Z$-equivalence. The book \cite{Newman} covers only the case of~$SL(2, \Z)$.
Our idea of classification is as follows. First, computations using programs CARAT (\cite{carat}, see also \cite{acta_cryst}) and GAP (\cite{gap}) allow to find the set of $16$ finite subgroups of $SL(3, \Z)$ containing representatives of all non-cyclic $\Z$-classes. They are listed in the table 1. (see theorem \ref{groups_table}). Next, in sections \ref{examples} and \ref{results}, we investigate their action on the product of three elliptic curves $A^3$ in order to compute cohomology of Kummer varieties. In this process we show that the actions defined by the listed groups are different. This section gives a short description of the computations and their results. The code of the programs can be found at www.mimuw.edu.pl/\~{}marysia/prog\_gap.

\subsection{General facts}
Let us first note some results which are base for understanding finite integral matrix groups. See e.g. \cite{Newman} and \cite{KP} for more general survey.

\begin{thm}[Minkowski]\label{Minkowski}
Let $G < GL(n, \Z)$ be a finite subgroup. Then, for every prime $p>2$ its $p$-th reduction $G \hookrightarrow GL(n, \Z) \rightarrow GL(n, \Z_p)$ is an embedding.
\end{thm}

\begin{lem}\label{existence_max_subgroups}
Every finite subgroup of $GL(n, \Z)$ is contained in a maximal finite subgroup.
\end{lem}

\begin{proof}
By Minkowski's theorem, the number of elements of a finite subgroup of $GL(n, \Z)$ divides the order of $GL(n, \Z_3)$. Hence an ascending sequence of finite subgroups must stabilize.
\end{proof}

By Minkowski's theorem, every finite subgroup of $GL(n, \Z)$ is isomorphic to a subgroup of $GL(n, \Z_3)$. Hence isomorphism types of finite subgroups of $SL(n, \Z)$ can be easily obtained for small values of $n$. The results for $SL(3, \Z)$ are mentioned in~\cite{Newman}, chapter~IX.14 or lemma~3.2 in~\cite{AW}.

\begin{lem}\label{types_of_subgroups}
Every nontrivial finite subgroup of $SL(3, \Z)$ is isomorphic to one of the following:
\begin{itemize}
\item cyclic groups $\Z_2$, $\Z_3$, $\Z_4$, $\Z_6$,
\item dihedral groups $D_4$, $D_6$, $D_8$, $D_{12}$,
\item the group $A_4$ of even permutations of four elements (e.g. the tetrahedral group~$T$),
\item the symmetric group $S_4$ of all permutations of four elements (e.g. the octahedral group~$O$).
\end{itemize}
\end{lem}

Note that $D_4$ is isomorphic to $\Z_2 \times \Z_2$, and $D_6$ to the symmetric group $S_3$.

The task is now to determine $\Z$-classes contained in each isomorphism class of finite subgroups of $SL(3, \Z)$. The algorithm is divided into two main steps, outlined in the next sections:
\begin{enumerate}
\item finding representatives of conjugacy classes of maximal finite subgroups of $GL(3, \Z)$,
\item determining $\Z$-classes of all finite subgroups of $SL(3, \Z)$ by analyzing the maximal subgroups of $GL(3, \Z)$.
\end{enumerate}

\subsection{Maximal finite subgroups in $GL(3, \Z)$}
To determine conjugacy classes of maximal finite subgroups of $GL(3, \Z)$ we used programs of the crystallographic package CARAT. Maximal finite subgroups of $GL(3, \Z)$ are so-called Bravais groups, which are integral matrix groups important in crystallography. This follows from the definition and basic properties of Bravais groups, described e.g. in \cite{acta_cryst} or \cite{BNZ}. CARAT programs can list $\Z$-classes of Bravais groups up to dimension~$6$. The only problem was to find maximal finite subgroups among them.

The program \emph{Bravais\_inclusions} (with option -S) listed all Bravais groups in $GL(3, \Z)$ containing $\{I, -I \} \izo \Z_2$. Every maximal subgroup of $GL(n, \Z)$ contains the matrix $-I$, because it belongs to the center of $GL(n, \Z)$, so all maximal finite subgroups of $GL(3, \Z)$ were on this list. The same program (with no option) let us cross out some groups which were not maximal elements of the list. There were four groups left: three of order $48$ and one of order $24$. It follows from the further part of the text that all of these are maximal finite subgroups of $GL(3, \Z)$ (their intersection with $SL(3, \Z)$ are not conjugate).

\subsection{All finite subgroups of $SL(3, \Z)$}
The following observation assures that we only need to consider subgroups contained in the intersections of maximal finite subgroups of $GL(3, \Z)$ with $SL(3, \Z)$. 

\begin{lem}
Let $G < SL(3, \Z)$ be a finite subgroup, and $\{H_i\}_{i\in I}$ is a set of representatives of all conjugacy classes of maximal finite subgroups of $GL(3, \Z)$. Then $G$ is $\Z$-equivalent to a $G' < H_i \cap SL(3, \Z)$ for some $i \in I$.
\end{lem}
\begin{proof}
By lemma \ref{existence_max_subgroups}, $G$ is contained in $H < GL(3, \Z)$ which is maximal, so conjugate to some $H_i$. Conjugation does not change the determinant of a matrix.
\end{proof}

We used the algebraic package GAP to determine $\Z$-classes of finite subgroups of $SL(3, \Z)$. 
In each maximal finite subgroup of $GL(3, \Z)$ its intersection with $SL(3, \Z)$ is of index $2$, so the input data consisted of three groups of order $24$ and one of order $12$. These are $S_4(1)$, $S_4(2)$, $S_4(3)$ and $D_{12}$ in table 1.

Our algorithm applied to each of these groups determines orders of their elements and created lists of subgroups of different isomorphism types by checking relations on generators. It can be implemented effectively, because all groups listed in the lemma \ref{types_of_subgroups} are cyclic or have presentation with two generators and one relation. As the first result we got $24$ groups but their number was decreased to~$20$ because in $4$ cases it was easy to find a base-change matrix.

Investigation of the action on $A^3$ suggested that $4$ of these groups were conjugate to some of the other $16$. Indeed, we found base-change matrices. The $16$ groups left represent different $\Z$-classes, which is proved by the analysis of their action on $A^3$ in Kummer construction described in sections \ref{examples} and \ref{results}. We note that Poincar\'e polynomials of Kummer varieties are not sufficient to distinguish $\Z$-classes of groups used in the construction. There are Kummer 3-folds obtained from non-conjugate groups, which have equal Poincar\'e polynomials (see table 2). In such cases we had to investigate the action on $A^3$ more carefully to distinguish $\Z$-classes. The most interesting cases \ref{case_D4_3} and \ref{case_D4_4} are discussed in section \ref{comment_D4}.

\begin{thm}\label{groups_table}
The list of groups in table 1 contains exactly one representative of each $\Z$-class of finite subgroups of $SL(3, \Z)$.
\end{thm}

\begin{table}[!hp]
\renewcommand{\arraystretch}{0.87}
\vspace{0.25cm}
\centerline{\begin{tabular}{|c|c|}
\hline
\textbf{group} & \textbf{generators}\\
\hline
$D_4(1)$ &
$\left(\begin{array}{rrr} -1 & 0 & 0 \\ 0 & -1 & 0 \\ 0 & 0 & 1 \end{array} \right)
\left(\begin{array}{rrr} 1 & 0 & 0 \\ 0 & -1 & 0 \\ 0 & 0 & -1 \end{array} \right)$\\
\hline
$D_4(2)$ &
$\left(\begin{array}{rrr} -1 & 0 & 0 \\ 0 & -1 & 0 \\ 0 & 0 & 1 \end{array} \right)
\left(\begin{array}{rrr} 0 & 1 & 0 \\ 1 & 0 & 0 \\ 0 & 0 & -1 \end{array} \right)$\\
\hline
$D_4(3)$ &
$\left(\begin{array}{rrr} -1 & 0 & 0 \\ 0 & 0 & -1 \\ 0 & -1 & 0 \end{array} \right)
\left(\begin{array}{rrr} 1 & 1 & 1 \\ 0 & -1 & 0 \\ 0 & 0 & -1 \end{array} \right)$\\
\hline
$D_4(4)$ &
$\left(\begin{array}{rrr} -1 & -1 & -1 \\ 0 & 0 & 1 \\ 0 & 1 & 0 \end{array} \right)
\left(\begin{array}{rrr} 0 & 0 & 1 \\ -1 & -1 & -1 \\ 1 & 0 & 0 \end{array} \right)$\\
\hline
$D_6(1)$ &
$\left(\begin{array}{rrr} -1 & 0 & 0 \\ 1 & 1 & 0 \\ 0 & 0 & -1 \end{array} \right)
\left(\begin{array}{rrr} -1 & -1 & 0 \\ 1 & 0 & 0 \\ 0 & 0 & 1 \end{array} \right)$\\
\hline
$D_6(2)$ &
$\left(\begin{array}{rrr} 0 & -1 & 0 \\ -1 & 0 & 0 \\ 0 & 0 & -1 \end{array} \right)
\left(\begin{array}{rrr} -1 & 1 & 0 \\ -1 & 0 & 0 \\ 0 & 0 & 1 \end{array} \right)$\\
\hline
$D_6(3)$ &
$\left(\begin{array}{rrr} -1 & 0 & 0 \\ 0 & 0 & -1 \\ 0 & -1 & 0 \end{array} \right)
\left(\begin{array}{rrr} 0 & -1 & 0 \\ 0 & 0 & 1 \\ -1 & 0 & 0 \end{array} \right)$\\
\hline
$D_8(1)$ &
$\left(\begin{array}{rrr} -1 & 0 & 0 \\ 0 & -1 & 0 \\ 0 & 0 & 1 \end{array} \right)
\left(\begin{array}{rrr} 0 & 0 & -1 \\ 0 & 1 & 0 \\ 1 & 0 & 0 \end{array} \right)$\\
\hline
$D_8(2)$ &
$\left(\begin{array}{rrr} -1 & -1 & -1 \\ 0 & 0 & 1 \\ 0 & 1 & 0 \end{array} \right)
\left(\begin{array}{rrr} 0 & -1 & 0 \\ 0 & 0 & -1 \\ 1 & 1 & 1 \end{array} \right)$\\
\hline
$D_{12}$ &
$\left(\begin{array}{rrr} 0 & 1 & 0 \\ 1 & 0 & 0 \\ 0 & 0 & -1 \end{array} \right)
\left(\begin{array}{rrr} 0 & 1 & 0 \\ -1 & 1 & 0 \\ 0 & 0 & 1 \end{array} \right)$\\
\hline
$A_4(1)$ &
$\left(\begin{array}{rrr} -1 & 0 & 0 \\ 0 & -1 & 0 \\ 0 & 0 & 1 \end{array} \right)
\left(\begin{array}{rrr} 0 & 0 & 1 \\ 1 & 0 & 0 \\ 0 & 1 & 0 \end{array} \right)$\\
\hline
$A_4(2)$ &
$\left(\begin{array}{rrr} -1 & -1 & -1 \\ 0 & 0 & 1 \\ 0 & 1 & 0 \end{array} \right)
\left(\begin{array}{rrr} 0 & 0 & 1 \\ 1 & 0 & 0 \\ 0 & 1 & 0 \end{array} \right)$\\
\hline
$A_4(3)$ &
$\left(\begin{array}{rrr} -1 & 0 & 0 \\ -1 & 0 & 1 \\ -1 & 1 & 0 \end{array} \right)
\left(\begin{array}{rrr} 0 & 0 & 1 \\ 1 & 0 & 0 \\ 0 & 1 & 0 \end{array} \right)$\\
\hline
$S_4(1)$ &
$\left(\begin{array}{rrr} 0 & 1 & 0 \\ 1 & 0 & 0 \\ 0 & 0 & -1 \end{array} \right)
\left(\begin{array}{rrr} 0 & 0 & 1 \\ 1 & 0 & 0 \\ 0 & 1 & 0 \end{array} \right)$\\
\hline
$S_4(2)$ &
$\left(\begin{array}{rrr} -1 & 0 & 0 \\ 0 & -1 & 0 \\ 1 & 1 & 1 \end{array} \right)
\left(\begin{array}{rrr} 0 & 0 & 1 \\ 1 & 0 & 0 \\ 0 & 1 & 0 \end{array} \right)$\\
\hline
$S_4(3)$ &
$\left(\begin{array}{rrr} -1 & 0 & 0 \\ 0 & 0 & -1 \\ 0 & -1 & 0 \end{array} \right)
\left(\begin{array}{rrr} -1 & 0 & 1 \\ -1 & 1 & 0 \\ -1 & 0 & 0 \end{array} \right)$\\
\hline
\end{tabular}}
\vspace{0.25cm}
\renewcommand{\arraystretch}{1}
\caption{$\Z$-classes of finite non-cyclic subgroups of $SL(3,\Z)$}
\end{table}

In the next sections, by abuse of notation, the names introduced in the table will denote both $\Z$-classes and chosen representatives.

\subsection{Relations of groups}\label{relations_groups}
The following definition is of considerable importance for the results of sections \ref{examples} and \ref{results}.

\begin{defn}\label{def_duality}
We say that finite subgroups $G, G' < GL(n, \Z)$ are dual, if by transposing all matrices in $G$ we obtain all matrices in $G'$. Conjugacy classes in $GL(n, \Z)$ are dual, if dual representatives can be chosen.
\end{defn}

\begin{prop}\label{dual_pairs}
Programs in GAP package allow to determine duality relation in the set of $\Z$-classes of finite subgroups of $SL(3, \Z)$:
\begin{itemize}
\item each of the following classes is dual to itself: $D_4(1)$, $D_4(2)$, $D_6(3)$, $D_8(1)$, $D_8(2)$, $D_{12}$, $A_4(1)$, $S_4(1)$;
\item there are four pairs of dual classes: $D_4(3)$ and $D_4(4)$, $D_6(1)$ and $D_6(2)$, $A_4(2)$ and $A_4(3)$, $S_4(2)$ and $S_4(3)$.
\end{itemize}
\end{prop}

Investigation of rational maps between Kummer 3-folds may lead to some new results. Therefore, by the following remark, we would like to understand the relation of inclusion (up to $\Z$-equivalence) of finite subgroups of $SL(3,\Z)$.

\begin{rem}\label{rational_map}
Let $G, H < GL(r, \Z)$ be finite subgroups such that there is $H' < G$ $\Z$-equivalent to H. Then there exists a rational map between Kummer varieties for $H$ and $G$.
\end{rem}

\begin{prop}\label{groups_relations_diagram}
The following diagram presents inclusions of finite non-cyclic subgroups of $SL(3, \Z)$ up to $\Z$-equivalence, determined by programs in GAP. An arrow from $H$ to $G$ means that there exists $H'<G$ which is $\Z$-equivalent to $H$. We omit arrows which come from composition of other arrows.

\vspace{3.25cm}
\vbox to 0ex{\vss\centerline{\includegraphics[width=0.6\textwidth]{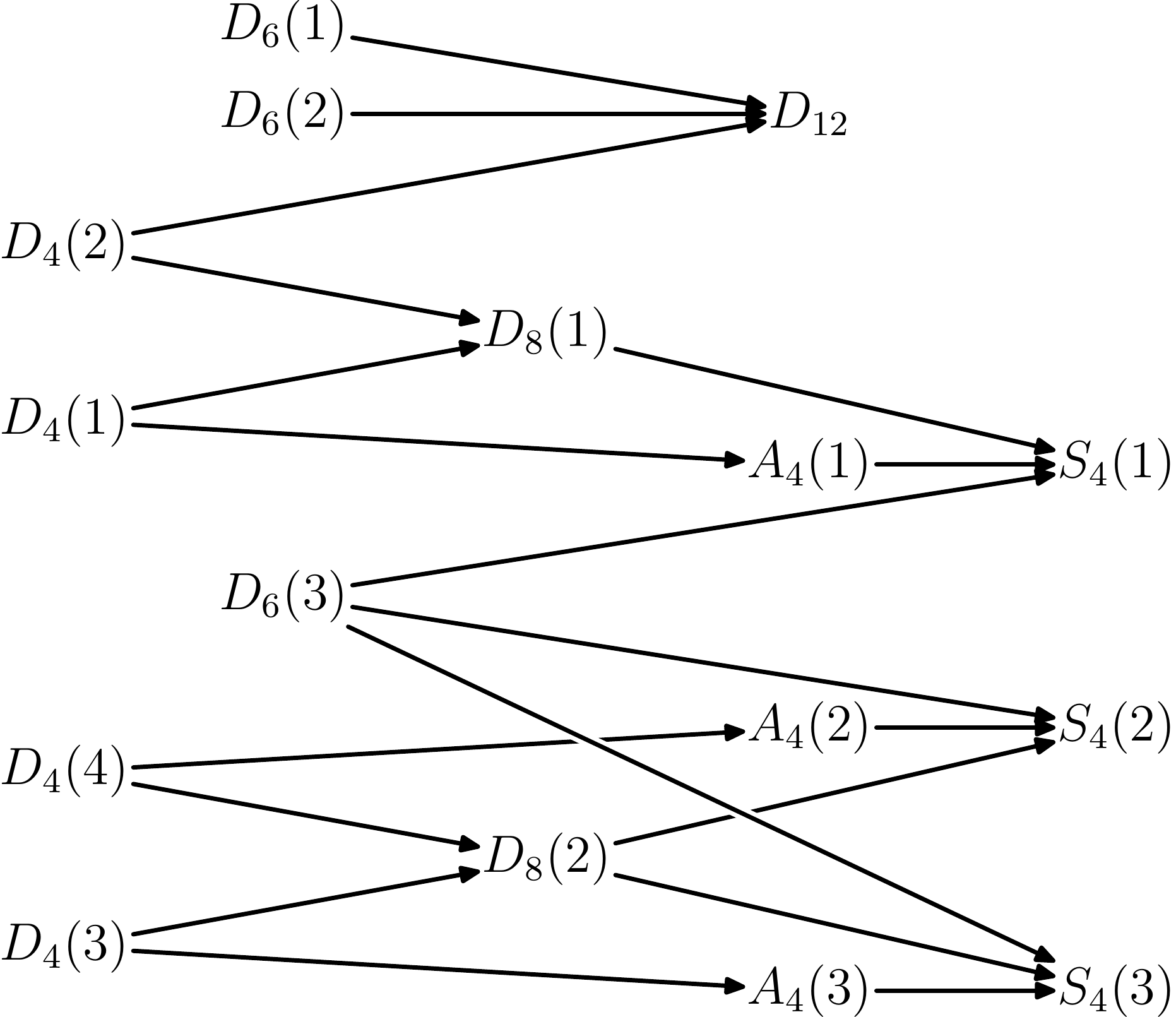}}\vss}
\vspace{3cm}
\end{prop}


\section{Construction and cohomology of Kummer 3-folds}\label{construction_details}

We compute Poincar\'e polynomial for the varieties obtained by applying Kummer construction to chosen representatives of all $\Z$-classes of finite subgroups of $SL(3, \Z)$. The general process of the construction and ideas of computations are described in \cite{AW}. This section is detailed discussion of the $3$-dimensional case. We describe the structure of Kummer 3-folds and explain methods of computing their cohomology, following the notation of \cite{AW}.

From now on, $G$ denotes a finite subgroup of $SL(3, \Z)$. For $H < G$ by $N(H)$ we denote the normalizer of $H$ in $G$ and $W(H) = N(H)/H$ is its Weyl group. To shorten the notation, if $H$ is cyclic, i.e. $H = \gen{h}$, we write $N(h)$ and $W(h)$ instead of $N(\gen{h})$ and $W(\gen{h})$.

\subsection{Stratification}\label{stratification}
We compute Poincar\'e polynomial $P_Y(t)$ of $Y = A^3/G$ and add the contribution of cohomology coming from the chosen resolution $f: X \ra Y$. By theorem~III.7.2 in~\cite{Bredon}, to determine Poincar\'e polynomial of $Y$ it is sufficient to compute dimensions of the spaces of $G$-invariant forms on $A^3$. This can be done by a simple function in GAP based on the character formulas given in \cite{FH}, chapter~2. Thus the main difficulty is understanding the contribution of the resolution of singularities.

We solve this problem using virtual Poincar\'e polynomials (see e.g. \cite{Arapura},~17.2 or \cite{Fulton},~4.5). Hence stratifications of $Y$ and $X$ must be chosen, as explained in \cite{AW},~2.2. There is a natural stratification of $Y = A^3/G$ determined by the action of~$G$ on $A^3$. The finite set of orbits of points with non-cyclic isotropy is the $0$-dimensional stratum. The $1$-dimensional stratum is the set of orbits of points with cyclic isotropy. The sum of $0$~and $1$-dimensional strata consist exactly of all singular points in $Y$, so the $3$-dimensional stratum is the smooth part of $Y$. By taking inverse images of the strata in $Y$ the decomposition of $X$ into $1$,~$2$ and $3$-dimensional strata is obtained. Note that the resolution of singularities $f: X \ra Y$ gives an isomorphism of $3$-dimensional strata.
The contribution to cohomology coming from the resolution is expressed by the difference of virtual Poincar\'e polynomials of the sums of $0$ and $1$-dimensional strata in~$Y$ and~$X$.

\subsection{Cohomology of the strata in $Y$}\label{cohomology_Y}
Let us first look at the set of points with non-trivial isotropy in the action of $G$ on $A^3$. By lemma 3.3 in \cite{AW}, the set $(A^3)^H$ of fixed points of a cyclic group $H$ consists of disjoint elliptic curves. The curves determined by $H$ can contain some points with non-cyclic isotropy, which are intersection points with curves determined by other cyclic groups. Hence the components of the set $(A^3)_0^H$ of points with isotropy $H$ are elliptic curves with some points removed. The image of $(A^3)_0^H$ in $Y$, denoted $Y([H])$, depends only on the conjugacy class of~$H$.

The Weyl group $W(H)$ acts freely on $(A^3)_0^H$. Let $K$ be a component of $Y([H])$ and $A_K$ an elliptic curve in $A^3$, which is mapped to $\ovl{K}$, the closure of $K$. By $W_K$ we denote the subgroup of $W(H)$ which fixes $A_K$, hence acts on the set of its points. Then the normalized closure of $K$, denoted $\widehat{K}$, is isomorphic to $A_K/W_K$ (see \cite{AW}, section 2.2). The virtual Poincar\'e polynomial of $K$ is the polynomial of $A_K/W_K$ with number of orbits of points on~$\ovl{K}$ with non-cyclic isotropy subtracted. There are three possible actions of $W_K$ on elliptic curve in our computations: trivial, the action of $\Z_2$ by involution, and the action of $\Z_2\times\Z_2$ generated by translation and involution. In the first case the Poincar\'e polynomial is $P_{A_K/W_K}(t) = 1 + 2t + t^2$. In the remaining two $A_K/W_K \izo \P^1$ and $P_{A_K/W_K}(t) = 1 + t^2$.

To obtain virtual Poincar\'e polynomial $P_K(t)$ of a component $K$ we have to subtract the number of points with non-cyclic isotropy on~$\ovl{K}$. Counting these points is generally easy, but there are two possible difficulties (which are in fact more important in analyzing the strata in $X$). The first is that sometimes the action of $W_K$ identifies some points with non-cyclic isotropy (see \ref{case_S4_2}). Then normalization of $\ovl{K}$ is not an isomorphism, because points with non-cyclic isotropy are not normal in $\ovl{K}$, but it is a homeomorphism. The second is when two points with non-cyclic isotropy lying on $A_K$ are mapped to the same point of the quotient, but the identification does not come from the action of $W_K$ (see \ref{case_S4_3}). Then $\ovl{K}$ goes twice through the image of such points. However, these situations appear only in some cases for $|G| \geq 8$, because they require $G$ having non-cyclic proper subgroup.

Summing $P_K(t)$ for all singular curves $K \subset Y$ we compute virtual Poincar\'e of $1$-dimensional stratum. Counting points in $0$-dimensional stratum is standard. It requires only the knowledge of non-cyclic subgroups of $G$ and their normalizers. Thus we compute the virtual Poincar\'e polynomial $P_3(t)$ of $3$-dimensional stratum both in $Y$ and in $X$, subtracting virtual Poincar\'e polynomials of $0$ and $1$-dimensional strata in $Y$ from $P_Y(t)$.

\subsection{Cohomology of the strata in $X$}\label{cohomology_X}
We now turn to computing the virtual Poincar\'e polynomial $P_2(t)$ of $2$-dimensional stratum in~$X$, which is the resolution of singularities over $Y([H])$ for all conjugacy classes of cyclic subgroups $H < G$. We consider locally product resolution $f : X \ra Y$ in the sense of definition~2.4 in~\cite{AW} (in fact, the resolution must be locally product in the case of Kummer $3$-folds, because the resolution is uniquely defined in codimension $2$). Let $F(H)$ be the fiber of minimal resolution of the quotient singularity $\C^2/H$, and $K$ a component of $Y[H]$. By lemma~2.5 in~\cite{AW}, to obtain the virtual Poincar\'e polynomial of $f^{-1}(K)$ we compute Poincar\'e polynomial of the quotient of $A_K\times F(H)$ by the action of $W_K$ and subtract Poincar\'e polynomial of the sum of fibers over points with non-cyclic isotropy.

Here we use the notation of $W_K$-Poincar\'e polynomials $P_{A_K\times F(H), W_K}(t)$ which have a representation of $W_K$ on $H^i(A_K\times F(H), \C)$ as a coefficient at $t^i$ (see \cite{AW}, section 2.1). To obtain Poincar\'e polynomial of the quotient $(A_K\times F(H))/W_K$ we apply the operation $\mu_0$ of taking the dimension of maximal trivial subrepresentations of coefficients to $P_{A_K\times F(H), W_K}(t)$. Hence 
$$P_{(A_K\times F(H))/W_K} = \mu_0(P_{A_K\times F(H), W_K}(t)) = \mu_0(P_{A_K, W_K}(t) \cdot P_{F(H), W_K}(t)),$$ 
which means that the task is now to compute $W_K$-Poincar\'e polynomials of~$A_K$ and~$F(H)$. 

The action of $W_K$ on $A_K$ is described in \ref{cohomology_Y}. The only remaining problem is the action of $W_K$ on $F(H)$. However, we can use the McKay correspondence for the minimal resolution of $\C^2/H$ (see \cite{Reid02}): the action of $W(H)$, hence also of $W_K$, on cohomology of $F(H)$ is the same as the action on conjugacy classes of $H$.

In our computations $H$ is one of the groups $\Z_2$, $\Z_3$, $\Z_4$, $\Z_6$. If $H = \Z_{n+1}$, the quotient singularity is of type $A_n$, so the non-zero cohomology spaces of $F(H)$ are $H^0(F(H)) = \C$ and $H^2(F(H)) = \C^n$. By McKay correspondence, the action of $W_K$ on a chosen basis of $H^2(F(H))$ is the same as the action of $W_K$ by conjugation on the set of nontrivial conjugacy classes in $H$. In fact we investigate the action of $W_K$ on $H$, because $H$ is abelian.

There are three cases: $W_K = 0$, $W_K = \Z_2$ and $W_K = \Z_2 \times \Z_2$. The last one appears three times in the computations, and only when $H = \Z_2$. Hence in the first and the last case $W_K$ acts trivially on $H^2(F(H))$. The same is true for~$W_K = \Z_2$ when~$H = \Z_2$.

Let us then look at the case $W_K = \Z_2$, assuming that $H \neq \Z_2$. Then $H$ has exactly two generators, which are interchanged by the action of $W_K$, because in our computations $H$ and $W_K$ always generate a non-abelian group. We choose new basis of $H^2(F(H))$: instead of each pair of basis vectors $\alpha$, $\beta$ interchanged by $W_K$ we take $\alpha + \beta$ and $\alpha - \beta$. In the new basis the representations of $W_K = \Z_2$ on $H^2(F(H))$ are as follows. The trivial representation is denoted by $1$, the standard by $\eps$ (that is $\eps = -1$) and the sum is direct sum of representations.

\vspace{0.25cm}
\centerline{\begin{tabular}{|c|c|c|}
\hline
\textbf{group $H$} & \textbf{representation} \\
\hline
$\Z_2$ & $1$ \\ \hline
$\Z_3$ & $1 + \eps$ \\ \hline
$\Z_4$ & $2 + \eps$ \\ \hline
$\Z_6$ & $3 + 2\eps$ \\ \hline
\end{tabular}}
\vspace{0.25cm}

The last step is computing the (virtual) Poincar\'e polynomial $P_1(t)$ of $1$-dimensional stratum of $X$, which is the resolution of quotient singularities $\C^3/H$ for a non-cyclic $H<G$. To understand this stratum it is sufficient to analyze a representative of each conjugacy class of non-cyclic subgroups of $G$. However, to compute $P_2(t)$ we need to count the points with non-cyclic singularity on each curve, which is much easier if we look at all non-cyclic subgroups of $G$, as in the following sections.

As for the existence of the crepant resolution of non-cyclic quotient singularities in $3$-dimensional Kummer construction, we rely on \cite{AW}, section 3.2. Again, to determine Poincar\'e polynomial we use McKay correspondence (this case is discussed in \cite{BM}): the number of $\P^1$ curves in the fiber of the resolution of quotient singularity $\C^3/H$ is equal to the number of nontrivial conjugacy classes in~$H$. This can be computed by a simple function in GAP. Results for all non-cyclic groups which appear in $3$-dimensional Kummer construction are summarized in the following table.

\vspace{0.25cm}
\centerline{\begin{tabular}{|c|c|}
\hline
\textbf{group $H$} & \textbf{Poincar\'e polynomial} \\
\hline
$D_4$ & $1 + 3t^2$ \\ \hline
$D_6$ & $1 + 2t^2$ \\ \hline
$D_8$ & $1 + 4t^2$ \\ \hline
$D_{12}$ & $1 + 5t^2$ \\ \hline
$A_4$ & $1 + 3t^2$ \\ \hline
$S_4$ & $1 + 4t^2$ \\ \hline
\end{tabular}}
\vspace{0.25cm}


\section{Examples}\label{examples}

Computations of Poincar\'e polynomials for all $16$ discussed groups are based on the methods described above. Therefore it is worth to look into the details only in the most complex cases, which contain all that is needed to be explained more carefully, because the remaining cases can be treated in the same way. In the present section we discuss Kummer construction for $S_4(2)$ and $S_4(3)$. The results of computations for the remaining groups are given in the next section.

\subsection*{Presentation of the results}
The information concerning curves of singular points is summarized in tables. In columns we give the following data:
\begin{itemize}
\item \textbf{group} --- isomorphism type of the investigated cyclic subgroup, (the type of its quotient singularity in brackets);
\item \textbf{generator} --- matrix which generates the subgroup in question (we chose representative of conjugacy classes of subgroups and a generator);
\item \textbf{equations} --- equations of the fixed points set in coordinates $(e_1, e_2, e_3)$ in $A^3$;
\item \textbf{components} --- number of elliptic curves in the fixed points set;
\item $W(g)$ --- the Weyl group;
\item \textbf{quotient} --- curves obtained as quotients of components of the fixed points set by the Weyl group action, with numbers of curves in each isomorphism class;
\item $W_K$ --- subgroup of $W(g)$ that acts on a single curve (given separately for each class of curves in the previous column).
\end{itemize}

\subsection{Case of $S_4(2)$}\label{case_S4_2}
There are three $\Z$-classes of $S_4$ subgroups in $SL(3, \Z)$. The case $S_4(1)$ (i.e. octahedral group) is treated in \cite{AW} and the action of this group on $A^3$ has a little simpler structure than those of the other two representations of $S_4$. Here we discuss $S_4(2)$, in some sense the most tricky one.

Let us recall some information about the structure of the group $S_4$. There are $2$ conjugacy classes of $\Z_2$ subgroups, one of all squares of order $4$ elements. All subgroups isomorphic to $\Z_4$ are conjugate, subgroups of type $\Z_3$ as well. Subgroups of type $D_4$ divide into two classes: one contains only a normal subgroup, with Weyl group $D_6$, and the second has three elements with Weyl group $\Z_2$. All $4$ subgroups isomorphic to $D_6$ are conjugate, as well as all $3$ subgroups of type $D_8$. There is also a normal subgroup $A_4$.

Next we choose a representative of the $\Z$-class $S_4(2)$: $$G = \left< \left(\begin{array}{rrr} -1 & 0 & 0 \\ 0 & -1 & 0 \\ 1 & 1 & 1 \end{array} \right),
\left(\begin{array}{rrr} 0 & 0 & 1 \\ 1 & 0 & 0 \\ 0 & 1 & 0 \end{array} \right)\right>.$$

The information about the singularities of $Y$ in codimension $2$ is easily obtained:

\vspace{0.25cm}
\centerline{\begin{tabular}{|c|c|c|c|c|c|c|}
\hline
\textbf{group} & \textbf{gen.} & \textbf{equ.} & \textbf{comp.} & $W(g)$ & \textbf{quot.} & $W_K$\\
\hline
$\Z_2$ $(A_1)$ &
$\left(\begin{array}{rrr} -1 & 0 & 0 \\ 0 & -1 & 0 \\ 1 & 1 & 1 \end{array} \right)$ &
$\begin{array}{c} 2e_1 = 0 \\ e_1 = e_2\\ \end{array}$ &
$4$ &
$\Z_2$ &
$4 \times \P^1$ &
$\Z_2$ \\
\hline
$\Z_2$ $(A_1)$ &
$\left(\begin{array}{rrr} 0 & 0 & 1 \\ -1 & -1 & -1 \\ 1 & 0 & 0 \end{array} \right)$ &
$\begin{array}{c} e_1 = e_3 \\ 2e_2 = -2e_1 \\ \end{array}$ &
$4$ &
$\Z_2\times \Z_2$ &
$3 \times \P^1$ &
$\Z_2 \times \Z_2$ \\
\hline
$\Z_3$ $(A_2)$ &
$\left(\begin{array}{rrr} 0 & 0 & 1 \\ 1 & 0 & 0 \\ 0 & 1 & 0 \end{array} \right)$ &
$\begin{array}{c} e_1 = e_3 \\ e_1 = e_2 \\ \end{array}$ &
$1$ &
$\Z_2$ &
$1 \times \P^1$ &
$\Z_2$ \\
\hline
$\Z_4$ $(A_3)$ &
$\left(\begin{array}{rrr} 0 & -1 & 0 \\ 0 & 0 & -1 \\ 1 & 1 & 1 \end{array} \right)$ &
$\begin{array}{c} e_1 = -e_2 \\ e_1 = e_3 \\ \end{array}$ &
$1$ &
$\Z_2$ &
$1 \times \P^1$ &
$\Z_2$ \\
\hline
\end{tabular}}
\vspace{0.25cm}

The most interesting part is analyzing the relations between $0$ and $1$ dimensional strata in $Y$ and finding the Poincar\'e polynomials for all strata after resolving the singularities.

It is much easier if we determine points with non-cyclic isotropy for all groups, not only for representatives of conjugacy classes.

Each of $3$ not normal $D_4$ subgroups fixes a set of $16$ points. The first set contain points of coordinates $(\alpha, \alpha, \beta)$, the second --- $(\alpha, \beta, \alpha)$, and the third --- $(\beta, \alpha, \alpha)$, where $2\alpha = 2\beta = 0$. Note that elements of the intersection, that is $4$ points $(\alpha, \alpha, \alpha)$, have isotropy $S_4$ and in each set the remaining $12$ points have isotropy $D_4$. The Weyl group $\Z_2$ acts on them, identifying pairs of points with $\alpha$ and $\beta$ interchanged. Because we look at conjugate subgroups, the action of $\Z_3 < S_4$ identifies triples with cyclicly permuted coordinates, one point taken from each family. Hence $36$ points with isotropy $D_4$ in $A^3$ are mapped to $6$ points of $Y$.

The normal $D_4$ subgroup fixes points which satisfy $e_1 = e_2 = e_3$ and $4e_1 = 0$, but these are exactly the points with isotropy $A_4$ ($12$ points with $2e_1 \neq 0$), or $S_4$. The set of points with isotropy $A_4$ is mapped to $6$ points of $Y$.

Points fixed by $D_6$ and $D_8$ subgroups satisfy $e_1 = e_2 = e_3$ and $2e_1 = 0$, so in fact they have isotropy $S_4$. Four points given by these equations are the only fixed points of $G$. Summing up, in the quotient there are $6$ points with quotient singularity of type $D_4$, $6$ with $A_4$ and $4$ with $S_4$.

We turn to the task of understanding resolutions over singular curves. Curves for the first (in the table) $\Z_2$ subgroup present what appears most often in the computations. Each component contains $3$ points with isotropy $D_4$ and $1$ with isotropy $S_4$. They are fixed by the action of $W_K$ (involution on each curve) and mapped to different points of $Y$. The fiber of resolution is $\P^1$, so $W_K$ acts trivially on it. We compute the virtual Poincar\'e polynomials using the notation in \ref{cohomology_X}. Note that from the polynomial of the quotient of a trivial bundle over elliptic curve we have to subtract the polynomial of the quotient of the induced bundle over the set of points with non-cyclic isotropy. The polynomial of a stratum of the resolution over one curve is
\begin{align*}
\mu_0((1 + \eps \cdot 2t + t^2)(1 + t^2)) - \mu_0(4(1+t^2)) = -3 -2t^2 + t^4.
\end{align*}

The curve of points with isotropy $\Z_4$ contains only $4$ points fixed by $S_4$, its virtual Poincar\'e polynomial is also easy to write:
\begin{align*}
\mu_0((1 + \eps \cdot 2t + t^2)(1 + (2 + \eps)t^2)) - \mu_0(4(1+(2 + \eps)t^2)) = -3 -5t^2 + 2t^3 + 2t^4.
\end{align*}

The second $\Z_2$ subgroup is one of three cases in our computations where $W_K = \Z_2 \times \Z_2$. These are isomorphic to $\P^1$, but the quotient map restricted to one component (after removing points with non-cyclic isotropy) is a $4$-sheeted cover. However, it is enough to look at the action of $W_K$ on the cohomology spaces of a single component. One generator of the product $\Z_2 \times \Z_2$ acts on elliptic curve $A$ by translation, which has trivial tangent map. The second, whichever we choose, acts by involution. The induced action on $H^0$ and $H^2$ is identity, and on $H^1$ has no fixed points, so the Poincar\'e polynomial is, indeed, $1 + t^2$. In all cases with $W_K = \Z_2 \times \Z_2$ the singularities are of type $A_1$, so we do not have to analyze the action of $\Z_2 \times \Z_2$ on the fibers.

Each of $3$ curves for the discussed subgroup contains $4$ points with isotropy $D_4$ and $4$ with isotropy $A_4$. Each set is mapped to $2$ points of the quotient. For the first set there is a stabilizer $\Z_2 < \Z_2 \times \Z_2$, so the representation of $W_K$ is $2(1 + \eps)$. The same is true for the second set, only the stabilizer is a different $\Z_2$ subgroup of $W_K$. Hence the Poincar\'e polynomial of one component is
\begin{align*}
\mu_0((1 + \eps \cdot 2t + t^2)(1 + t^2)) - \mu_0(4(1 + \eps)(1+t^2)) = -3 - 2t^2 + t^4.
\end{align*}

And finally the only case where it can be non-obvious to compute correctly the Poincar\'e polynomial. The main difficulty is computing the polynomial of the quotient of the bundle over the set of points with non-cyclic isotropy. In this single case the representations of $W_K$ both on this set of points and on the fiber are nontrivial, which has to be taken into account. The curve of fixed points of $\Z_3$ contains $12$ points with isotropy $A_4$, which are mapped to $6$ points of the quotient, so the representation of $W_K = \Z_2$ on this set is $6(1 + \eps)$ and the $W_K$-Poincar\'e polynomial is $6(1 + \eps)(1 + (1+\eps)t^2)$.
There are also $4$ fixed points of $S_4$. The Poincar\'e polynomial~is

\begin{align*}
\mu_0((1 + \eps \cdot 2t + t^2)(1 + (1 + \eps)t^2)) - \mu_0((6(1 + \eps) + 4)(1+(1 + \eps)t^2)) =\\= t^4 + 2t^3 - 14t^2 - 9.
\end{align*}

For $S_4(2)$, and also for other investigated representations of $S_4$, the (virtual) Poincar\'e polynomial of the quotient $Y$ is $$P_Y(t) = t^6 + t^4 + 4t^3 + t^2 + 1.$$
The polynomial of $3$-dimensional stratum can be easily computed using $P_Y$ and the given data. The polynomial of the $2$-dimensional stratum is a sum of the polynomials for all components. And the polynomials of the resolutions over points with non-cyclic isotropy, which sum up to $1$-dimensional stratum, can be taken from the list in \ref{cohomology_X}.
Hence the virtual Poincar\'e polynomials of the strata are the following:
\begin{align*}
P_3(t) &= P_Y(t) - 8(1 + t^2 -4) - (1 + t^2 - 10) - 16 = t^6 + t^4 + 4t^3 - 8t^2 + 18,\\
P_2(t) &= 4(1 + t^2 - 4)(1 + t^2) + 3\mu_0((1 + t^2 - 4(1 + \eps))(1 + t^2)) +\\ &\quad + \mu_0((1 + \eps \cdot 2t + t^2 - 4 - 6(1 + \eps))(1 + (1 + \eps)t^2)) + \\ &\quad + \mu_0((1 + \eps \cdot 2t + t^2 - 4)(1 + (2 + \eps)t^2)) = 10t^4 + 4t^3 - 33t^2 - 33, \\
P_1(t) &= 12(1 + 3t^2) + 4(1 + 4t^2) = 52t^2 + 16,
\end{align*}
and, finally, the Poincar\'e polynomial of $X$ is $$P_X(t) = t^6 + 11t^4 + 8t^3 + 11t^2 + 1.$$

\subsection{Case of $S_4(3)$}\label{case_S4_3}
There is one more interesting detail which does not appear in the previous example, but can be observed in the case of $S_4(3)$. Hence we explain it more carefully, while on the other steps of computations only general information is given.

We choose $G = \left< \left(\begin{array}{rrr} -1 & 0 & 0 \\ 0 & 0 & -1 \\ 0 & -1 & 0 \end{array} \right),
\left(\begin{array}{rrr} -1 & 0 & 1 \\ -1 & 1 & 0 \\ -1 & 0 & 0 \end{array} \right)\right>$.

\vspace{0.25cm}
\centerline{\begin{tabular}{|c|c|c|c|c|c|c|}
\hline
\textbf{group} & \textbf{gen.} & \textbf{equ.} & \textbf{comp.} & $W(g)$ & \textbf{quot.} & $W_K$\\
\hline
$\Z_2$ $(A_1)$ &
$\left(\begin{array}{rrr} -1 & 0 & 0 \\ 0 & 0 & -1 \\ 0 & -1 & 0 \end{array} \right)$ &
$\begin{array}{c} 2e_1 = 0 \\ e_2 = -e_3\\ \end{array}$ &
$4$ &
$\Z_2$ &
$4 \times \P^1$ &
$\Z_2$ \\
\hline
$\Z_2$ $(A_1)$ &
$\left(\begin{array}{rrr} -1 & 0 & 0 \\ -1 & 0 & 1 \\ -1 & 1 & 0 \end{array} \right)$ &
$\begin{array}{c} 2e_1 = 0 \\ e_3 = e_1 + e_2 \\ \end{array}$ &
$4$ &
$\Z_2\times \Z_2$ &
$3 \times \P^1$ &
$\Z_2 \times \Z_2$ \\
\hline
$\Z_3$ $(A_2)$ &
$\left(\begin{array}{rrr} -1 & 0 & 1 \\ -1 & 1 & 0 \\ -1 & 0 & 0 \end{array} \right)$ &
$\begin{array}{c} e_1 = 0 \\ e_3 = 0 \\ \end{array}$ &
$1$ &
$\Z_2$ &
$1 \times \P^1$ &
$\Z_2$ \\
\hline
$\Z_4$ $(A_3)$ &
$\left(\begin{array}{rrr} 0 & -1 & 1 \\ 0 & 0 & 1 \\ -1 & 0 & 1 \end{array} \right)$ &
$\begin{array}{c} e_1 = 0 \\ e_2 = e_3 \\ \end{array}$ &
$1$ &
$\Z_2$ &
$1 \times \P^1$ &
$\Z_2$ \\
\hline
\end{tabular}}
\vspace{0.25cm}

In $A^3$ there are $3$ sets of $12$ points with isotropy $D_4$, mapped to $6$ points of the quotient, two points from each set to one in $Y$. There are also $6$ points determined by the normal $D_4$ subgroup, identified in the quotient. The points with isotropy $D_6$ are divided into $4$ sets of $3$ points, mapped to $3$ points of $Y$, because all subgroups of type $D_6$ are conjugate. Similarly there are $3$ sets of $3$ points with isotropy $D_8$. The action of normalizer is trivial, sets are identified by the quotient map. One point has isotropy $S_4$.

In the case of $S_4(3)$ there are curves which go twice through some points. Let us look at the components of the fixed points set of the first $\Z_2$ subgroup. Three of them contain $4$ points with isotropy $D_4$ each. Take one component and denote it $A_K$, as in \ref{cohomology_Y} and \ref{cohomology_X}. On $A_K$ points with isotropy $D_4$ are fixed by $W_K$, but they are mapped to $2$ points of the quotient. That is, the identification does not come from the action of $W_K$, so the neighborhoods of these points are not glued by the quotient map. It means that if we take the image of $A_K$ in $Y$ and remove from it $2$ points with singularity $D_4$, we get $K$ isomorphic to $\P^1$ without $4$ points. In other words, the quotient of an elliptic curve $A_K$ by $W_K = \Z_2$ is not isomorphic to the closure of $K$ in $Y$, but only to its normalized closure $\widehat{K}$.

Apart from $4$ points with isotropy $D_4$, three of the curves for the first $\Z_2$ subgroup contain also $2$ points with isotropy $D_6$ each. They are identified by the action of $W_K$. The fourth curve contains $3$ points with isotropy $D_8$ and $1$ with $S_4$, mapped to different points of $Y$.

The quotient map on components for the second $\Z_2$, with $W_K = \Z_2 \times \Z_2$, is $4$-sheeted cover, as in the previous example. Each of these curves contains $6$ points with isotropy $D_4$ and $2$ with $D_8$, which are pairwise identified by the action of $W_K$. The curve for $\Z_3$ subgroup contains $3$ points with isotropy $D_6$, mapped to different points of the quotient, and $1$ fixed point of $S_4$. The curve for $\Z_4$ contains $3$ points with isotropy $D_8$, different in $Y$, and $1$ fixed by $S_4$.

Using the same methods as in the case of $S_4(2)$ we give the virtual Poincar\'e polynomials of all strata in $Y$ with resolved singularities:
\begin{align*}
P_3(t) &= P_Y(t) - 3(1 + t^2 - 5) - 6(1 + t^2 - 4) - 14 = t^6 + t^4 + 4t^3 - 8t^2 + 17,\\
P_2(t) &= (1 + t^2 - 4)(1 + t^2) + 3\mu_0((1 + \eps \cdot 2t + t^2 - 4 - (1 + \eps))(1 + t^2)) +\\ &\quad + 
3\mu_0((1 + \eps \cdot 2t + t^2 - 4(1 + \eps))(1 + t^2)) +\\ &\quad + \mu_0((1 + \eps \cdot 2t + t^2 - 4)(1 + (1 + \eps)t^2)) +\\ &\quad +
\mu_0((1 + \eps \cdot 2t + t^2 - 4)(1 + (2 + \eps)t^2)) = 10t^4 + 4t^3 - 24t^2 - 30,\\
P_1(t) &= 7(1 + 3t^2) + 3(1 + 2t^2) + 4(1 + 4t^2) = 43t^2 + 14.
\end{align*}
The Poincar\'e polynomial of $X$ is $$P_X(t) = t^6 + 11t^4 + 8t^3 + 11t^2 + 1.$$


\section{Results of the computations}\label{results}
In this section we collect the results of computations for all Kummer varieties except the examples from the previous section. Moreover, we provide some information about the action of finite subgroups of $SL(3, \Z)$ on $A^3$, structure of their quotients (curves of singular points, their equations) and virtual Poincar\'e polynomials for all strata.

\subsection{Cases of $D_4$}
For all of the investigated representations of $D_4$ the (virtual) Poincar\'e polynomial of the quotient $Y$ is $$P_Y(t) = t^6 + 3t^4 + 8t^3 + 3t^2 + 1.$$

The group $D_4$ has three normal subgroups of order $2$, which determine curves of points of the isotropy $\Z_2$. The isotropy of their intersection points is $D_4$.

\subsubsection{$D_4(1)$}\label{case_D4_1}
\enlargethispage{0.5\baselineskip}
$G = \left<\left(\begin{array}{rrr} -1 & 0 & 0 \\ 0 & -1 & 0 \\ 0 & 0 & 1 \end{array} \right),
\left(\begin{array}{rrr} 1 & 0 & 0 \\ 0 & -1 & 0 \\ 0 & 0 & -1 \end{array} \right)\right>$

\vspace{0.25cm}
\centerline{\begin{tabular}{|c|c|c|c|c|c|c|}
\hline
\textbf{group} & \textbf{gen.} & \textbf{equ.} & \textbf{comp.} & $W(g)$ & \textbf{quot.} & $W_K$\\
\hline
$\Z_2$ $(A_1)$ &
$\left(\begin{array}{rrr} -1 & 0 & 0 \\  0 & -1 & 0 \\ 0 & 0 & 1 \end{array}\right)$ &
$\begin{array}{c} 2e_1 = 0 \\ 2e_2 = 0\\ \end{array}$ &
$16$ &
$\Z_2$ &
$16\times \P^1$ &
$\Z_2$ \\
\hline
$\Z_2$ $(A_1)$ &
$\left(\begin{array}{rrr} 1 & 0 & 0 \\  0 & -1 & 0 \\ 0 & 0 & -1 \end{array}\right)$ &
$\begin{array}{c} 2e_2 = 0 \\ 2e_3 = 0\\ \end{array}$ &
$16$ &
$\Z_2$ &
$16\times \P^1$ &
$\Z_2$ \\
\hline
$\Z_2$ $(A_1)$ &
$\left(\begin{array}{rrr} -1 & 0 & 0 \\  0 & 1 & 0 \\ 0 & 0 & -1 \end{array}\right)$ &
$\begin{array}{c} 2e_1 = 0 \\ 2e_3 = 0\\ \end{array}$ &
$16$ &
$\Z_2$ &
$16\times \P^1$ &
$\Z_2$ \\
\hline
\end{tabular}}
\vspace{0.25cm}

There are $64$ points with isotropy $D_4$, $4$ on each curve. They are defined by equations $2e_1 = 2e_2 = 2e_3 = 0$.
The virtual Poincar\'e polynomials of the strata are
\begin{align*}
P_3(t) &= P_Y(t) - 48(1 + t^2 - 4) - 64 = t^6 + 3t^4 + 8t^3 - 45t^2 + 81,\\
P_2(t) &= 48(1 + t^2 -4)(1 + t^2) = 48t^4 - 96t^2 - 144,\\
P_1(t) &= 64(1 + 3t^2) = 192t^2 + 64,
\end{align*}
and the Poincar\'e polynomial of $X$ is $$P_X(t) = t^6 + 51t^4 + 8t^3 + 51t^2 + 1.$$

\subsubsection{$D_4(2)$}\label{case_D4_2}
$G = \left<\left(\begin{array}{rrr} -1 & 0 & 0 \\ 0 & -1 & 0 \\ 0 & 0 & 1 \end{array} \right),
\left(\begin{array}{rrr} 0 & 1 & 0 \\ 1 & 0 & 0 \\ 0 & 0 & -1 \end{array} \right)\right>$

\vspace{0.25cm}
\centerline{\begin{tabular}{|c|c|c|c|c|c|c|}
\hline
\textbf{group} & \textbf{gen.} & \textbf{equ.} & \textbf{comp.} & $W(g)$ & \textbf{quot.} & $W_K$\\
\hline
$\Z_2$ $(A_1)$ &
$\left(\begin{array}{rrr} -1 & 0 & 0 \\ 0 & -1 & 0 \\ 0 & 0 & 1  \end{array}\right)$ &
$\begin{array}{c} 2e_1 = 0 \\ 2e_2 = 0\\ \end{array}$ &
$16$ &
$\Z_2$ &
$\begin{array}{c} 4 \times \P^1 \\ 6 \times A\\ \end{array}$ &
$\begin{array}{c} \Z_2 \\ 0\\ \end{array}$ \\
\hline
$\Z_2$ $(A_1)$ &
$\left(\begin{array}{rrr} 0 & 1 & 0 \\ 1 & 0 & 0 \\ 0 & 0 & -1 \end{array}\right)$ &
$\begin{array}{c} e_1 = e_2 \\ 2e_3 = 0\\ \end{array}$ &
$4$ &
$\Z_2$ &
$4\times \P^1$ &
$\Z_2$ \\
\hline
$\Z_2$ $(A_1)$ &
$\left(\begin{array}{rrr} 0 & -1 & 0 \\ -1 & 0 & 0 \\ 0 & 0 & -1 \end{array}\right)$ &
$\begin{array}{c} e_1 = -e_2 \\ 2e_3 = 0\\ \end{array}$ &
$4$ &
$\Z_2$ &
$4\times \P^1$ &
$\Z_2$ \\
\hline
\end{tabular}}
\vspace{0.25cm}

There are $16$ points with isotropy $D_4$. They are given by equations $2e_1 = 2e_2 = 2e_3 = 0$ and $e_1 = e_2$. Each of the $\P^1$ curves contains $4$ of these points. The elliptic curves do not contain any of them.
The virtual Poincar\'e polynomials of the strata are
\begin{align*}
P_3(t) &= P_Y(t) - 12(1 + t^2 - 4) - 6(1 + 2t + t^2) - 16 =\\ &= t^6 + 3t^4 + 8t^3 - 15t^2 -12t + 15,\\
P_2(t) &= 12(1 + t^2 -4)(1 + t^2) + 6(1 + 2t + t^2)(1 + t^2) =\\ &= 18t^4 + 12t^3 - 12t^2 + 12t - 30,\\
P_1(t) &= 16(1 + 3t^2) = 48t^2 + 16,
\end{align*}
and the Poincar\'e polynomial of $X$ is $$P_X(t) = t^6 + 21t^4 + 20t^3 + 21t^2 + 1.$$

\subsubsection{$D_4(3)$}\label{case_D4_3}
$G = \left<\left(\begin{array}{rrr} -1 & 0 & 0 \\ 0 & 0 & -1 \\ 0 & -1 & 0 \end{array} \right)
\left(\begin{array}{rrr} 1 & 1 & 1 \\ 0 & -1 & 0 \\ 0 & 0 & -1 \end{array} \right)\right>$

\vspace{0.25cm}
\centerline{\begin{tabular}{|c|c|c|c|c|c|c|}
\hline
\textbf{group} & \textbf{gen.} & \textbf{equ.} & \textbf{comp.} & $W(g)$ & \textbf{quot.} & $W_K$\\
\hline
$\Z_2$ $(A_1)$ &
$\left(\begin{array}{rrr} -1 & 0 & 0 \\ 0 & 0 & -1 \\ 0 & -1 & 0  \end{array}\right)$ &
$\begin{array}{c} 2e_1 = 0 \\ e_2 = -e_3\\ \end{array}$ &
$4$ &
$\Z_2$ &
$4 \times \P^1$ &
$\Z_2$ \\
\hline
$\Z_2$ $(A_1)$ &
$\left(\begin{array}{rrr} 1 & 1 & 1 \\ 0 & -1 & 0 \\ 0 & 0 & -1 \end{array}\right)$ &
$\begin{array}{c} 2e_2 = 0 \\ e_2 = e_3\\ \end{array}$ &
$4$ &
$\Z_2$ &
$4\times \P^1$ &
$\Z_2$ \\
\hline
$\Z_2$ $(A_1)$ &
$\left(\begin{array}{rrr} -1 & -1 & -1 \\ 0 & 0 & 1 \\ 0 & 1 & 0 \end{array}\right)$ &
$\begin{array}{c} e_2 = e_3 \\ 2e_1 = -2e_2\\ \end{array}$ &
$4$ &
$\Z_2$ &
$4\times \P^1$ &
$\Z_2$ \\
\hline
\end{tabular}}
\vspace{0.25cm}

There are $16$ points with isotropy $D_4$, $4$ on each curve. They are given by equations $2e_1 = 2e_2 = 0$ and $e_2 = e_3$.
The virtual Poincar\'e polynomials of the strata are
\begin{align*}
P_3(t) &= P_Y(t) - 12(1 + t^2 - 4) - 16 = t^6 + 3t^4 + 8t^3 - 9t^2 + 21,\\
P_2(t) &= 12(1 + t^2 -4)(1 + t^2) = 12t^4 - 24t^2 - 36,\\
P_1(t) &= 16(1 + 3t^2) = 48t^2 + 16,
\end{align*}
and the Poincar\'e polynomial of $X$ is $$P_X(t) = t^6 + 15t^4 + 8t^3 + 15t^2 + 1.$$

\subsubsection{$D_4(4)$}\label{case_D4_4}
$G = \left<\left(\begin{array}{rrr} -1 & -1 & -1 \\ 0 & 0 & 1 \\ 0 & 1 & 0 \end{array} \right)
\left(\begin{array}{rrr} 0 & 0 & 1 \\ -1 & -1 & -1 \\ 1 & 0 & 0 \end{array} \right)\right>$

\vspace{0.25cm}
\centerline{\begin{tabular}{|c|c|c|c|c|c|c|}
\hline
\textbf{group} & \textbf{gen.} & \textbf{equ.} & \textbf{comp.} & $W(g)$ & \textbf{quot.} & $W_K$\\
\hline
$\Z_2$ $(A_1)$ &
$\left(\begin{array}{rrr} -1 & -1 & -1 \\ 0 & 0 & 1 \\ 0 & 1 & 0 \end{array}\right)$ &
$\begin{array}{c} e_2 = e_3 \\ 2e_1 = -2e_2\\ \end{array}$ &
$4$ &
$\Z_2$ &
$4 \times \P^1$ &
$\Z_2$ \\
\hline
$\Z_2$ $(A_1)$ &
$\left(\begin{array}{rrr} 0 & 0 & 1 \\ -1 & -1 & -1 \\ 1 & 0 & 0 \end{array}\right)$ &
$\begin{array}{c} e_1 = e_3 \\ 2e_2 = -2e_1\\ \end{array}$ &
$4$ &
$\Z_2$ &
$4\times \P^1$ &
$\Z_2$ \\
\hline
$\Z_2$ $(A_1)$ &
$\left(\begin{array}{rrr} 0 & 1 & 0 \\ 1 & 0 & 0 \\ -1 & -1 & -1 \end{array}\right)$ &
$\begin{array}{c} e_1 = e_2 \\ 2e_3 = -2e_1\\ \end{array}$ &
$4$ &
$\Z_2$ &
$4\times \P^1$ &
$\Z_2$ \\
\hline
\end{tabular}}
\vspace{0.25cm}

There are $16$ points with isotropy $D_4$, $4$ on each curve. They are defined by equations $e_1 = e_2 = e_3$ and $4e_1 = 0$.
The virtual Poincar\'e polynomials of the strata are
\begin{align*}
P_3(t) &= P_Y(t) - 12(1 + t^2 - 4) - 16 = t^6 + 3t^4 + 8t^3 - 9t^2 + 21,\\
P_2(t) &= 12(1 + t^2 -4)(1 + t^2) = 12t^4 - 24t^2 - 36,\\
P_1(t) &= 16(1 + 3t^2) = 48t^2 + 16,
\end{align*}
and the Poincar\'e polynomial of $X$ is $$P_X(t) = t^6 + 15t^4 + 8t^3 + 15t^2 + 1.$$

\subsubsection{How do the cases \ref{case_D4_3} and \ref{case_D4_4} differ?}\label{comment_D4}
In the cases \ref{case_D4_3} and \ref{case_D4_4} not only Poincar\'e polynomials, but also numbers of singular curves and points are the same. However, these are different cases, because in the case of $D_4(3)$ the set of points with nontrivial isotropy in $X$ is connected, which is not true for $D_4(4)$. In this section $x$ stands for the arbitrary point of elliptic curve $A$ and the letters $\alpha$ and $\beta$ denote points which satisfy $2x=0$.

In the case of $D_4(3)$, points with isotropy $D_4$ are $(\alpha, \beta, \beta)$ for all possible values of $\alpha$ and $\beta$. Components of fixed points set for the first $\Z_2$ subgroup are parametrized by $(\alpha, x, -x)$. For the second $\Z_2$ parametrizations of components are $(x, \beta, \beta)$, so each of these curves intersect all curves determined by the first $\Z_2$. Therefore the sum of $0$ and $1$ dimensional strata is connected.

In the case of $D_4(4)$ there are three families of the fixed points curves for $\Z_2$ subgroups, each containing $4$ curves, one for each value of $\alpha$. They are parametrized by $(\alpha - x, x, x)$, $(x, \alpha - x, x)$ and $(x, x, \alpha - x)$ respectively. Points with isotropy $D_4$ are $(\gamma, \gamma, \gamma)$, where $4\gamma = 0$. A choice of $\alpha$ determines three curves, one from each family, which intersect in $4$ points such that $2\gamma = \alpha$. Curves for different values $\alpha_0$ and $\alpha_1$ do not intersect. Hence the sum of $0$ and $1$ dimensional strata has $4$ connected components, which shows that the groups $D_4(3)$ and $D_4(4)$ cannot be conjugate.

\begin{figure}[h]
\vspace{2.75cm}
\vbox to 0ex{\vss\centerline{\includegraphics[width=0.75\textwidth]{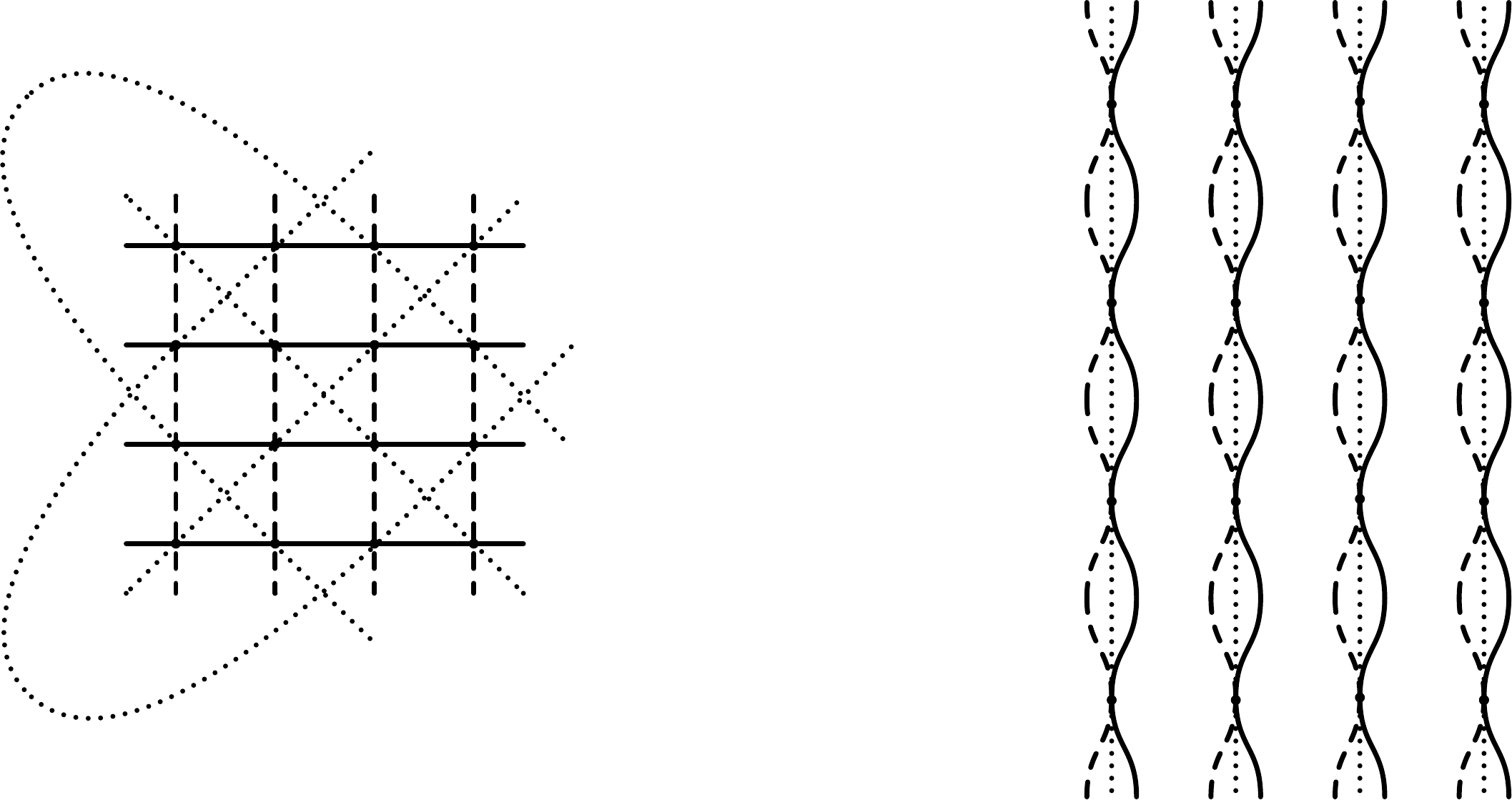}}\vss} 
\vspace{3cm}
\caption{Structure of the sets of singular points in $Y$ for $D_4(3)$ and $D_4(4)$}
\end{figure}

Figure 1. presents sums of $0$ and $1$ dimensional strata in discussed cases. Intersection points of curves are marked by black dots. Components determined by the same $\Z_2$ subgroup are drawn in the same line style.

Note that by transposing matrices in $D_4(3)$ we obtain a group conjugate to $D_4(4)$.

\subsection{Cases of $D_6$}
For all of the investigated representations of $D_6 \izo S_3$ the (virtual) Poincar\'e polynomial of the quotient $Y$ is
$$P_Y(t) = t^6 + 2t^4 + 6t^3 + 2t^2 + 1.$$

The symmetric group $S_3$ has a normal subgroup of order $3$. Elements of order $2$ determine three conjugate subgroups $\Z_2$. Points which have non-cyclic isotropy are fixed points of the action of $S_3$.

\subsubsection{$D_6(1)$}\label{case_D6_1}
$G = \left< \left(\begin{array}{rrr} -1 & 0 & 0 \\ 1 & 1 & 0 \\ 0 & 0 & -1 \end{array} \right),
\left(\begin{array}{rrr} -1 & -1 & 0 \\ 1 & 0 & 0 \\ 0 & 0 & 1 \end{array} \right)\right>$

\vspace{0.25cm}
\centerline{\begin{tabular}{|c|c|c|c|c|c|c|}
\hline
\textbf{group} & \textbf{gen.} & \textbf{equ.} & \textbf{comp.} & $W(g)$ & \textbf{quot.} & $W_K$\\
\hline
$\Z_2$ $(A_1)$ &
$\left(\begin{array}{rrr} -1 & 0 & 0 \\ 1 & 1 & 0 \\ 0 & 0 & -1 \end{array} \right)$ &
$\begin{array}{c} e_1 = 0 \\ 2e_3 = 0\\ \end{array}$ &
$4$ &
$0$ &
$4 \times A$ &
$0$ \\
\hline
$\Z_3$ $(A_2)$ &
$\left(\begin{array}{rrr} -1 & -1 & 0 \\ 1 & 0 & 0 \\ 0 & 0 & 1 \end{array} \right)$ &
$\begin{array}{c} e_1 = e_2 \\ 3e_1 = 0 \\ \end{array}$ &
$9$ &
$\Z_2$ &
$\begin{array}{c} 4\times A \\ 1\times \P^1 \\ \end{array}$ &
$\begin{array}{c} 0 \\ \Z_2 \\ \end{array}$ \\
\hline
\end{tabular}}
\vspace{0.25cm}

Fixed points of the action of $D_6$ are given by equations $e_1 = e_2 = 0$ and $2e_3 = 0$, so there are $4$ of them. All lie on the $\P^1$ curve of fixed points of $\Z_3$. Each elliptic curve for $\Z_2$ goes through exactly one of these points. Elliptic curves for $\Z_3$ do not contain any points of bigger isotropy.
The virtual Poincar\'e polynomials of the strata are
\begin{align*}
P_3(t) &= P_Y(t) - 4(1 + 2t + t^2 - 1) - 4(1 + 2t + t^2) - (1 + t^2 - 4) - 4 =\\ &= t^6 + 2t^4 + 6t^3 - 7t^2 - 16t - 4,\\
P_2(t) &= 4(1 + 2t + t^2 - 1)(1 + t^2) + 4(1 + 2t + t^2)(1 + 2t^2) +\\ &\quad + \mu_0((1 + \eps \cdot 2t + t^2 - 4)(1 + (1 + \eps)t^2)) = 
13t^4 + 26t^3 + 14t^2 + 16t + 1,\\
P_1(t) &= 4(1 + 2t^2) = 8t^2 + 4,
\end{align*}
and the Poincar\'e polynomial of $X$ is $$P_X(t) = t^6 + 15t^4 + 32t^3 + 15t^2 + 1.$$

\subsubsection{$D_6(2)$}\label{case_D6_2}
$G = \left< \left(\begin{array}{rrr} 0 & -1 & 0 \\ -1 & 0 & 0 \\ 0 & 0 & -1 \end{array} \right),
\left(\begin{array}{rrr} -1 & 1 & 0 \\ -1 & 0 & 0 \\ 0 & 0 & 1 \end{array} \right)\right>$

\vspace{0.25cm}
\centerline{\begin{tabular}{|c|c|c|c|c|c|c|}
\hline
\textbf{group} & \textbf{gen.} & \textbf{equ.} & \textbf{comp.} & $W(g)$ & \textbf{quot.} & $W_K$\\
\hline
$\Z_2$ $(A_1)$ &
$\left(\begin{array}{rrr} 0 & -1 & 0 \\ -1 & 0 & 0 \\ 0 & 0 & -1 \end{array} \right)$ &
$\begin{array}{c} e_1 = -e_2 \\ 2e_3 = 0\\ \end{array}$ &
$4$ &
$0$ &
$4 \times A$ &
$0$ \\
\hline
$\Z_3$ $(A_2)$ &
$\left(\begin{array}{rrr} -1 & 1 & 0 \\ -1 & 0 & 0 \\ 0 & 0 & 1 \end{array} \right)$ &
$\begin{array}{c} 2e_1 = e_2 \\ 3e_1 = 0 \\ \end{array}$ &
$9$ &
$\Z_2$ &
$9 \times \P^1$ &
$\Z_2$ \\
\hline
\end{tabular}}
\vspace{0.25cm}

Fixed points of the action of $D_6$ are given by equations $3e_1 = 0$, $e_2 = 2e_1$ and $2e_3 = 0$, so there are $36$ of them. There are $9$ of them on each elliptic curve for $\Z_2$ and $4$ on each $\P^1$ for $\Z_3$.
The virtual Poincar\'e polynomials of the strata are the following:
\begin{align*}
P_3(t) &= P_Y(t) - 4(1 + 2t + t^2 - 9) - 9(1 + t^2 - 4) - 36 =\\ &= t^6 + 2t^4 + 6t^3 - 11t^2 - 8t + 24,\\
P_2(t) &= 4(1 + 2t + t^2 - 9)(1 + t^2) + 9\mu_0((1 + \eps \cdot 2t + t^2 - 4)(1 + (1 + \eps)t^2)) = \\ &=
13t^4 + 26t^3 - 46t^2 + 8t - 59,\\
P_1(t) &= 36(1 + 2t^2) = 72t^2 + 36,
\end{align*}
and the Poincar\'e polynomial of $X$ is $$P_X(t) = t^6 + 15t^4 + 32t^3 + 15t^2 + 1.$$

\subsubsection{$D_6(3)$}\label{case_D6_3}
$G = \left< \left(\begin{array}{rrr} -1 & 0 & 0 \\ 0 & 0 & -1 \\ 0 & -1 & 0 \end{array} \right),
\left(\begin{array}{rrr} 0 & -1 & 0 \\ 0 & 0 & 1 \\ -1 & 0 & 0 \end{array} \right)\right>$

\vspace{0.25cm}
\centerline{\begin{tabular}{|c|c|c|c|c|c|c|}
\hline
\textbf{group} & \textbf{gen.} & \textbf{equ.} & \textbf{comp.} & $W(g)$ & \textbf{quot.} & $W_K$\\
\hline
$\Z_2$ $(A_1)$ &
$\left(\begin{array}{rrr} -1 & 0 & 0 \\ 0 & 0 & -1 \\ 0 & -1 & 0 \end{array} \right)$ &
$\begin{array}{c} 2e_1 = 0 \\ e_2 = -e_3\\ \end{array}$ &
$4$ &
$0$ &
$4 \times A$ &
$0$ \\
\hline
$\Z_3$ $(A_2)$ &
$\left(\begin{array}{rrr} 0 & -1 & 0 \\ 0 & 0 & 1 \\ -1 & 0 & 0 \end{array} \right)$ &
$\begin{array}{c} e_1 = -e_2 \\ e_2 = e_3 \\ \end{array}$ &
$1$ &
$\Z_2$ &
$1\times \P^1$ &
$\Z_2$ \\
\hline
\end{tabular}}
\vspace{0.25cm}

Fixed points of the action are defined by equations $2e_1 = 0$ and $e_1 = e_2 = e_3$, so there are $4$ of them. Each lies on the $\P^1$ curve of fixed points of the $\Z_3$ action and on one of the elliptic curves for $\Z_2$.
The virtual Poincar\'e polynomials of the strata are the following:
\begin{align*}
P_3(t) &= P_Y(t) - 4(1 + 2t + t^2 - 1) - (1 + t^2 - 4) - 4 = t^6 + 2t^4 + 6t^2 - 3t^2 - 8t,\\
P_2(t) &= 4(1 + 2t + t^2 - 1)(1 + t^2) + \mu_0((1 + \eps \cdot 2t + t^2 - 4)(1 + (1 + \eps)t^2)) = \\ &= 5t^4 + 10t^3 + 2t^2 + 8t - 3,\\
P_1(t) &= 4(1 + 2t^2) = 8t^2 + 4,
\end{align*}
and the Poincar\'e polynomial of $X$ is $$P_X(t) = t^6 + 7t^4 + 16t^3 + 7t^2 + 1.$$

\subsection{Cases of $D_8$}
For both investigated representations of $D_6 \izo S_3$ the (virtual) Poincar\'e polynomial of the quotient $Y$ is
$$P_Y(t) = t^6 + 2t^4 + 6t^3 + 2t^2 + 1.$$

The group $D_8$ has normal subgroup $\Z_4$, which contains a $\Z_2$ subgroup, also normal in $D_8$, with Weyl group $\Z_2\times\Z_2$. There are two other classes of $\Z_2$ subgroups, each containing $2$ groups, with Weyl group $\Z_2$. Points with non-cyclic isotropy are fixed by the whole $D_8$ or only by one of $2$ non-conjugate subgroups isomorphic to $D_4$.

\subsubsection{$D_8(1)$}\label{case_D8_1}
$G = \left< \left(\begin{array}{rrr} -1 & 0 & 0 \\ 0 & -1 & 0 \\ 0 & 0 & 1 \end{array} \right),
\left(\begin{array}{rrr} 0 & 0 & -1 \\ 0 & 1 & 0 \\ 1 & 0 & 0 \end{array} \right)\right>$\\

\vspace{0.25cm}
\centerline{\begin{tabular}{|c|c|c|c|c|c|c|}
\hline
\textbf{group} & \textbf{gen.} & \textbf{equ.} & \textbf{comp.} & $W(g)$ & \textbf{quot.} & $W_K$\\
\hline
$\Z_2$ $(A_1)$ &
$\left(\begin{array}{rrr} -1 & 0 & 0 \\ 0 & -1 & 0 \\ 0 & 0 & 1 \end{array} \right)$ &
$\begin{array}{c} 2e_1 = 0 \\ 2e_2 = 0\\ \end{array}$ &
$16$ &
$\Z_2$ &
$16 \times \P^1$ &
$\Z_2$ \\
\hline
$\Z_2$ $(A_1)$ &
$\left(\begin{array}{rrr} 0 & 0 & -1 \\ 0 & -1 & 0 \\ -1 & 0 & 0 \end{array} \right)$ &
$\begin{array}{c} e_1 = -e_3 \\ 2e_2 = 0\\ \end{array}$ &
$4$ &
$\Z_2$ &
$4 \times \P^1$ &
$\Z_2$ \\
\hline
$\Z_2$ $(A_1)$ &
$\left(\begin{array}{rrr} -1 & 0 & 0 \\ 0 & 1 & 0 \\ 0 & 0 & -1 \end{array} \right)$ &
$\begin{array}{c} 2e_1 = 0 \\ 2e_3 = 0\\ \end{array}$ &
$16$ &
$\Z_2\times\Z_2$ &
$6 \times \P^1$ &
$\Z_2$ \\
\hline
$\Z_4$ $(A_3)$ &
$\left(\begin{array}{rrr} 0 & 0 & -1 \\ 0 & 1 & 0 \\ 1 & 0 & 0 \end{array} \right)$ &
$\begin{array}{c} 2e_1 = 0 \\ e_1 = e_3 \\ \end{array}$ &
$4$ &
$\Z_2$ &
$4\times \P^1$ &
$\Z_2$ \\
\hline
\end{tabular}}
\vspace{0.25cm}

In this case $64$ points given by equations $2e_1 = 2e_2 = 2e_3 = 0$ have non-cyclic isotropy. There are $48$ points with isotropy $D_4$, which are mapped to $24$ points of the quotient, two to one. The remaining $16$ points are fixed by $D_8$.

Each $\P^1$ curve determined by the first $\Z_2$ subgroup contains $3$ points with isotropy $D_4$ and $1$ fixed by $D_8$. There are $4$ points with isotropy $D_8$ on each curve for the second $\Z_2$. Fixed point set of the normal $\Z_2$ subgroup has $16$ components. Four of them consists of points with isotropy $\Z_4$, we look now at the other $12$. The Weyl group $\Z_2 \times \Z_2$ acts on this set by involution on each component and identification of pairs of curves. Hence they are mapped to $6$ copies of $\P^1$ in $Y$ and the quotient map is double cover on each component. Each of $12$ curves contains $4$ points with isotropy $D_4$. The curves determined by $\Z_4$ contain $4$ points with isotropy $D_8$ each.

The virtual Poincar\'e polynomials of the strata are the following:
\begin{align*}
P_3(t) &= P_Y(t) - 30(1 + t^2 - 4) - 40 = t^6 + 2t^4 + 6t^3 - 28t^2 + 51,\\
P_2(t) &= 26(1 + t^2 - 4)(1 + t^2) + 4\mu_0((1 + \eps \cdot 2t + t^2 - 4)(1 + (2 + \eps)t^2)) =\\ &= 34t^4 + 8t^3 - 72t^2 - 90,\\
P_1(t) &= 24(1 + 3t^2) + 16(1 + 4t^2) = 136t^2 + 40,
\end{align*}
and the Poincar\'e polynomial of $X$ is $$P_X(t) = t^6 + 36t^4 + 14t^3 + 36t^2 + 1.$$

\subsubsection{$D_8(2)$}\label{case_D8_2}
$G = \left< \left(\begin{array}{rrr} -1 & -1 & -1 \\ 0 & 0 & 1 \\ 0 & 1 & 0 \end{array} \right),
\left(\begin{array}{rrr} 0 & -1 & 0 \\ 0 & 0 & -1 \\ 1 & 1 & 1 \end{array} \right)\right>$

\vspace{0.25cm}
\centerline{\begin{tabular}{|c|c|c|c|c|c|c|}
\hline
\textbf{group} & \textbf{gen.} & \textbf{equ.} & \textbf{comp.} & $W(g)$ & \textbf{quot.} & $W_K$\\
\hline
$\Z_2$ $(A_1)$ &
$\left(\begin{array}{rrr} -1 & -1 & -1 \\ 0 & 0 & 1 \\ 0 & 1 & 0 \end{array} \right)$ &
$\begin{array}{c} e_2 = e_3 \\ 2e_1 = -2e_2\\ \end{array}$ &
$4$ &
$\Z_2$ &
$4 \times \P^1$ &
$\Z_2$ \\
\hline
$\Z_2$ $(A_1)$ &
$\left(\begin{array}{rrr} -1 & 0 & 0 \\ 1 & 1 & 1 \\ 0 & 0 & -1 \end{array} \right)$ &
$\begin{array}{c} 2e_1 = 0 \\ e_3 = e_1\\ \end{array}$ &
$4$ &
$\Z_2$ &
$4 \times \P^1$ &
$\Z_2$ \\
\hline
$\Z_2$ $(A_1)$ &
$\left(\begin{array}{rrr} 0 & 0 & 1 \\ -1 & -1 & -1 \\ 1 & 0 & 0 \end{array} \right)$ &
$\begin{array}{c} e_1 = e_3 \\ 2e_2 = -2e_1\\ \end{array}$ &
$4$ &
$\Z_2\times\Z_2$ &
$3 \times \P^1$ &
$\Z_2 \times \Z_2$ \\
\hline
$\Z_4$ $(A_3)$ &
$\left(\begin{array}{rrr} 0 & -1 & 0 \\ 0 & 0 & -1 \\ 1 & 1 & 1 \end{array} \right)$ &
$\begin{array}{c} e_1 = -e_2 \\ e_1 = e_3 \\ \end{array}$ &
$1$ &
$\Z_2$ &
$1\times \P^1$ &
$\Z_2$ \\
\hline
\end{tabular}}
\vspace{0.25cm}

Each of $D_4$ subgroups fixes $16$ points. The first set is given by equations $e_1 = e_2 = e_3$ and $4e_1 = 0$, the second by $e_1 = e_3$ and $2e_1 = 2e_2 = 0$. The intersection of these sets consists of $4$ points with isotropy $D_8$. The remaining $24$ points with isotropy $D_4$ are mapped to $12$ points of the quotient, two to one.

Three of the curves determined by the first $\Z_2$ subgroup contain $4$ points with isotropy $D_4$ each. These points are pairwise identified in the quotient, but the identification does not come from the action of $W_K$ (as in \ref{case_S4_3}). The fourth curve contains $4$ points with isotropy $D_8$. Each curve for the second $\Z_2$ contains $3$ points with isotropy $D_4$ and one with isotropy $D_8$. The quotient map on the curves determined by the normal $\Z_2$ subgroup is a $4$-sheeted cover (as in \ref{case_S4_2}). There are $8$ points with isotropy $D_4$ on each of these curves, $4$ for each $D_4$ subgroup. Pairs of these points are identified by the action of a $\Z_2$ subgroup of Weyl group, as in the cases of $S_4(2)$ and $S_4(3)$. All points fixed by $D_8$ lie on the $\P^1$ determined by $\Z_4$.

The virtual Poincar\'e polynomials of the strata are the following:
\begin{align*}
P_3(t) &= P_Y(t) - 12(1 + t^2 - 4) - 16 = t^6 + 2t^4 + 6t^3 - 10t^2 + 21,\\
P_2(t) &= 11(1 + t^2 - 4)(1 + t^2) + \mu_0((1 + \eps \cdot 2t + t^2 - 4)(1 + (2 + \eps)t^2)) =\\ &= 13t^4 + 2t^3 - 27t^2 - 36,\\
P_1(t) &= 12(1 + 3t^2) + 4(1 + 4t^2) = 52t^2 + 16,
\end{align*}
and the Poincar\'e polynomial of $X$ is $$P_X(t) = t^6 + 15t^4 + 8t^3 + 15t^2 + 1.$$

\subsection{Case of $D_{12}$}\label{case_D12}

There is only one $\Z$-class of subgroups of $SL(3, \Z)$ isomorphic to $D_{12}$. The (virtual) Poincar\'e polynomial of the quotient $Y$ is
$$P_Y(t) = t^6 + 2t^4 + 6t^3 + 2t^2 + 1.$$

The group $D_{12}$ has normal subgroup $\Z_6$, which contains $\Z_3$ and $\Z_2$, also normal. Other $\Z_2$ subgroups divide into $2$ classes. There are also $3$ conjugate subgroups isomorphic to $D_4$ and $2$ normal $D_6$ subgroups.

The chosen representative is
$G = \left< \left(\begin{array}{rrr} 0 & 1 & 0 \\ 1 & 0 & 0 \\ 0 & 0 & -1 \end{array} \right),
\left(\begin{array}{rrr} 0 & 1 & 0 \\ -1 & 1 & 0 \\ 0 & 0 & 1 \end{array} \right)\right>$.

\vspace{0.25cm}
\centerline{\begin{tabular}{|c|c|c|c|c|c|c|}
\hline
\textbf{group} & \textbf{gen.} & \textbf{equ.} & \textbf{comp.} & $W(g)$ & \textbf{quot.} & $W_K$\\
\hline
$\Z_2$ $(A_1)$ &
$\left(\begin{array}{rrr} -1 & 0 & 0 \\ -1 & 1 & 0 \\ 0 & 0 & -1 \end{array} \right)$ &
$\begin{array}{c} e_1 = 0 \\ 2e_3 = 0\\ \end{array}$ &
$4$ &
$\Z_2$ &
$4 \times \P^1$ &
$\Z_2$ \\
\hline
$\Z_2$ $(A_1)$ &
$\left(\begin{array}{rrr} -1 & 1 & 0 \\ 0 & 1 & 0 \\ 0 & 0 & -1 \end{array} \right)$ &
$\begin{array}{c} 2e_1 = e_2 \\ 2e_3 = 0 \\ \end{array}$ &
$4$ &
$\Z_2$ &
$4 \times \P^1$ &
$\Z_2$ \\
\hline
$\Z_2$ $(A_1)$ &
$\left(\begin{array}{rrr} -1 & 0 & 0 \\ 0 & -1 & 0 \\ 0 & 0 & 1 \end{array} \right)$ &
$\begin{array}{c} 2e_1 = 0 \\ 2e_2 = 0 \\ \end{array}$ &
$16$ &
$D_6$ &
$\begin{array}{c} 3 \times \P^1 \\ 1 \times A \\ \end{array}$ &
$\begin{array}{c} \Z_2 \\ 0 \\ \end{array}$ \\
\hline
$\Z_3$ $(A_2)$ &
$\left(\begin{array}{rrr} -1 & 1 & 0 \\ -1 & 0 & 0 \\ 0 & 0 & 1 \end{array} \right)$ &
$\begin{array}{c} 2e_1 = e_2 \\ 3e_1 = 0 \\ \end{array}$ &
$9$ &
$\Z_2\times \Z_2$ &
$4 \times \P^1$ &
$\Z_2$ \\
\hline
$\Z_6$ $(A_5)$ &
$\left(\begin{array}{rrr} 0 & 1 & 0 \\ -1 & 1 & 0 \\ 0 & 0 & 1 \end{array} \right)$ &
$\begin{array}{c} e_1 = 0 \\ e_2 = 0 \\ \end{array}$ &
$1$ &
$\Z_2$ &
$1 \times \P^1$ &
$\Z_2$ \\
\hline
\end{tabular}}
\vspace{0.25cm}

There are $36$ points with isotropy $D_4$, mapped to $12$ points of the quotient (triples are identified by the action of $\Z_3$). Next, there are $32$ points with isotropy $D_6$, mapped to $16$ points of $Y$ (pairs are identified), and $4$ points fixed by $D_{12}$.

Each curve determined by the first $\Z_2$ subgroup contains $3$ points with isotropy $D_4$ and $1$ with isotropy $D_{12}$. Each curve for the second $\Z_2$ subgroup also contains $3$ points with isotropy $D_4$ and $1$ fixed by $D_{12}$, and moreover $8$ with isotropy $D_6$, which are mapped to $4$ points of $Y$ (identification comes from the action of $W_K$). As for the normal $\Z_2$ subgroup, there are $4$ points with isotropy $D_4$ on each of its $\P^1$ curves and no points with non-cyclic isotropy on its elliptic curve. The $\P^1$ curves determined by $\Z_3$ contain $4$ points with isotropy $D_6$ each. Finally, all points fixed by $D_{12}$ lie on the curve of fixed points of $\Z_6$.

The virtual Poincar\'e polynomials of the strata are the following:
\begin{align*}
P_3(t) &= P_Y(t) - 4(1 + t^2 - 4) - 4(1 + t^2 - 8) - 3(1 + t^2 - 4) - (1 + 2t + t^2) -\\ &\quad -4(1 + t^2 - 4) - (1 + t^2 - 4) - 32 = t^6 + 2t^4 + 6t^3 - 15t^2 - 2t + 32,\\
P_2(t) &= 7(1 - 4 + (2 - 4)t^2 + t^4) + 4(1 - 4 - 4 + (2 - 4 - 4)t^2 + t^4) +\\ &\quad + (1 + 2t + 2t^2 + 2t^3 + t^4) + 4(1 - 4 + (2 - 4)t^2 + 2t^3 + t^4) + \\ &\quad + (1 - 4 + ((1 - 4)\cdot 3 + 1)t^2 + 4t^3 + 3t^4) = 19t^4 + 14t^3 - 52t^2 + 2t - 63,\\
P_1(t) &= 12(1 + 3t^2) + 16(1 + 2t^2) + 4(1 + 5t^2) = 88t^2 + 32.
\end{align*}
and the Poincar\'e polynomial of $X$ is $$P_X(t) = t^6 + 21t^4 + 20t^3 + 21t^2 + 1.$$

\subsection{Cases of $A_4$}
For all of the investigated representations of $A_4$ the (virtual) Poincar\'e polynomial of the quotient $Y$ is
$$P_Y(t) = t^6 + t^4 + 4t^3 + t^2 + 1.$$

In $A_4$ all cyclic subgroups of order $3$ are conjugate, the same applies to the subgroups of order $2$. There is also one non-cyclic subgroup, containing all elements of order $4$, isomorphic to $D_4$ and normal.

\subsubsection{$A_4(1)$}\label{case_A4_1}
$G = \left< \left(\begin{array}{rrr} -1 & 0 & 0 \\ 0 & -1 & 0 \\ 0 & 0 & 1 \end{array} \right),
\left(\begin{array}{rrr} 0 & 0 & 1 \\ 1 & 0 & 0 \\ 0 & 1 & 0 \end{array} \right)\right>$

\vspace{0.25cm}
\centerline{\begin{tabular}{|c|c|c|c|c|c|c|}
\hline
\textbf{group} & \textbf{gen.} & \textbf{equ.} & \textbf{comp.} & $W(g)$ & \textbf{quot.} & $W_K$\\
\hline
$\Z_2$ $(A_1)$ &
$\left(\begin{array}{rrr} -1 & 0 & 0 \\ 0 & -1 & 0 \\ 0 & 0 & 1 \end{array} \right)$ &
$\begin{array}{c} 2e_1 = 0 \\ 2e_2 = 0\\ \end{array}$ &
$16$ &
$\Z_2$ &
$16 \times \P^1$ &
$\Z_2$ \\
\hline
$\Z_3$ $(A_2)$ &
$\left(\begin{array}{rrr} 0 & 0 & 1 \\ 1 & 0 & 0 \\ 0 & 1 & 0 \end{array} \right)$ &
$\begin{array}{c} e_1 = e_2 \\ e_2 = e_3 \\ \end{array}$ &
$1$ &
$0$ &
$1\times A$ &
$0$ \\
\hline
\end{tabular}}
\vspace{0.25cm}

There are $64$ points with isotropy $D_4$, given by equations $2e_1 = 2e_2 = 2e_3 = 0$. Elements of order $3$ act on this set by cyclic permutation of coordinates, so in the quotient triples of points are identified, except $4$ points fixed by the action of $A_4$. Hence the $0$-dimensional stratum consists of $4$ points with isotropy $A_4$, which lie on the elliptic curve of fixed points of $\Z_3$, and $20$ points with isotropy $D_4$.  Each curve of fixed points of $\Z_2$ contains $4$ points of $0$-dimensional stratum.

The virtual Poincar\'e polynomials of the strata are the following:
\begin{align*}
P_3(t) &= P_Y(t) - 16(1 + t^2 - 4) - (1 + 2t + t^2 - 4) - 24 =\\ &= t^6 + t^4 + 4t^3 - 16t^2 - 2t + 28,\\
P_2(t) &= 16(1 + t^2 - 4)(1 + t^2) + (1 + 2t + t^2 - 4)(1 + 2t^2) =\\ &= 18t^4 + 4t^3 -37t^2 + 2t - 51,\\
P_1(t) &= 20(1 + 3t^2) + 4(1 + 3t^2) = 72t^2 + 24,
\end{align*}
and the Poincar\'e polynomial of $X$ is $$P_X(t) = t^6 + 19t^4 + 8t^3 + 19t^2 + 1.$$

\subsubsection{$A_4(2)$}\label{case_A4_2}

$G = \left< \left(\begin{array}{rrr} -1 & -1 & -1 \\ 0 & 0 & 1 \\ 0 & 1 & 0 \end{array} \right),
\left(\begin{array}{rrr} 0 & 0 & 1 \\ 1 & 0 & 0 \\ 0 & 1 & 0 \end{array} \right)\right>$

\vspace{0.25cm}
\centerline{\begin{tabular}{|c|c|c|c|c|c|c|}
\hline
\textbf{group} & \textbf{gen.} & \textbf{equ.} & \textbf{comp.} & $W(g)$ & \textbf{quot.} & $W_K$\\
\hline
$\Z_2$ $(A_1)$ &
$\left(\begin{array}{rrr} -1 & -1 & -1 \\ 0 & 0 & 1 \\ 0 & 1 & 0 \end{array} \right)$ &
$\begin{array}{c} e_2 = e_3 \\ 2e_1 = -2e_2\\ \end{array}$ &
$4$ &
$\Z_2$ &
$4 \times \P^1$ &
$\Z_2$ \\
\hline
$\Z_3$ $(A_2)$ &
$\left(\begin{array}{rrr} 0 & 0 & 1 \\ 1 & 0 & 0 \\ 0 & 1 & 0 \end{array} \right)$ &
$\begin{array}{c} e_1 = e_2 \\ e_2 = e_3 \\ \end{array}$ &
$1$ &
$0$ &
$1\times A$ &
$0$ \\
\hline
\end{tabular}}
\vspace{0.25cm}

There are $16$ fixed points of $D_4$, they satisfy $e_1 = e_2 = e_3$ and $4e_1 = 0$. The chosen generator of $\Z_3$ acts trivially on this set, so these points are fixed by $A_4$. They all lie on the elliptic curve for $\Z_3$, and on each $\P^1$ of fixed points of $\Z_2$ there are $4$ of them.

The virtual Poincar\'e polynomials of the strata are the following:
\begin{align*}
P_3(t) &= P_Y(t) - 4(1 + t^2 - 4) - (1 + 2t + t^2 - 16) - 16 =\\ &= t^6 + t^4 + 4t^3 - 4t^2 - 2t + 12,\\
P_2(t) &= 4(1 + t^2 - 4)(1 + t^2) + (1 + 2t + t^2 - 16)(1 + 2t^2) =\\ &= 6t^4 + 4t^3 -37t^2 + 2t - 27,\\
P_1(t) &= 16(1 + 3t^2) = 48t^2 + 16,
\end{align*}
and the Poincar\'e polynomial of $X$ is $$P_X(t) = t^6 + 7t^4 + 8t^3 + 7t^2 + 1.$$

\subsubsection{$A_4(3)$}\label{case_A4_3}
$G = \left< \left(\begin{array}{rrr} -1 & 0 & 0 \\ -1 & 0 & 1 \\ -1 & 1 & 0 \end{array} \right),
\left(\begin{array}{rrr} 0 & 0 & 1 \\ 1 & 0 & 0 \\ 0 & 1 & 0 \end{array} \right)\right>$

\vspace{0.25cm}
\centerline{\begin{tabular}{|c|c|c|c|c|c|c|}
\hline
\textbf{group} & \textbf{gen.} & \textbf{equ.} & \textbf{comp.} & $W(g)$ & \textbf{quot.} & $W_K$\\
\hline
$\Z_2$ $(A_1)$ &
$\left(\begin{array}{rrr} -1 & 0 & 0 \\ -1 & 0 & 1 \\ -1 & 1 & 0 \end{array} \right)$ &
$\begin{array}{c} 2e_1 = 0 \\ e_3 = e_1 + e_2\\ \end{array}$ &
$4$ &
$\Z_2$ &
$4 \times \P^1$ &
$\Z_2$ \\
\hline
$\Z_3$ $(A_2)$ &
$\left(\begin{array}{rrr} 0 & 0 & 1 \\ 1 & 0 & 0 \\ 0 & 1 & 0 \end{array} \right)$ &
$\begin{array}{c} e_1 = e_2 \\ e_2 = e_3 \\ \end{array}$ &
$1$ &
$0$ &
$1\times A$ &
$0$ \\
\hline
\end{tabular}}
\vspace{0.25cm}

Fixed points of $D_4$ are described by equations $2e_1 = 2e_2 = 0$ and $e_3 = e_1 + e_2$; there are $16$ of them. One is fixed by $A_4$. The remaining $15$ are permuted by the action of $\Z_3$, so in the quotient we get $5$ points with the isotropy $D_4$. The elliptic curve of fixed points of $\Z_3$ contains the fixed point of $A_4$. On each $P^1$ curves there are $4$ points with non-cyclic isotropy. One of them contains the fixed point of $A_4$ and $3$ more points. Three of them contain $4$ points with isotropy $D_4$, two being identified in the quotient, but not as a result of the normalizer's action. Hence each of the images of these curves contain three points with non-cyclic isotropy, and goes two times through one of them.

The virtual Poincar\'e polynomials of the strata are the following:
\begin{align*}
P_3(t) &= P_Y(t) - 4(1 + t^2 - 4) - (1 + 2t + t^2 -1) - 6 = t^6 + t^4 + 4t^3 - 4t^2 - 2t + 7,\\
P_2(t) &= 4(1 + t^2 - 4)(1 + t^2) + (1 + 2t + t^2 - 1)(1 + 2t^2) =\\ &= 6t^4 + 4t^3 - 7t^2 + 2t -12,\\
P_1(t) &= 5(1 + 3t^2) + (1 + 3t^2) = 18t^2 + 6,
\end{align*}
and the Poincar\'e polynomial of $X$ is $$P_X(t) = t^6 + 7t^4 + 8t^3 + 7t^2 + 1.$$

\subsection{Case of $S_4(1)$}\label{case_S4_1}
$G = \left< \left(\begin{array}{rrr} 0 & 1 & 0 \\ 1 & 0 & 0 \\ 0 & 0 & -1 \end{array} \right),
\left(\begin{array}{rrr} 0 & 0 & 1 \\ 1 & 0 & 0 \\ 0 & 1 & 0 \end{array} \right)\right>$

This case is computed in \cite{AW}, but we repeat the results to complete the survey.

\vspace{0.25cm}
\centerline{\begin{tabular}{|c|c|c|c|c|c|c|}
\hline
\textbf{group} & \textbf{gen.} & \textbf{equ.} & \textbf{comp.} & $W(g)$ & \textbf{quot.} & $W_K$\\
\hline
$\Z_2$ $(A_1)$ &
$\left(\begin{array}{rrr} 0 & 1 & 0 \\ 1 & 0 & 0 \\ 0 & 0 & -1 \end{array} \right)$ &
$\begin{array}{c} e_1 = e_2 \\ 2e_3 = 0\\ \end{array}$ &
$4$ &
$\Z_2$ &
$4 \times \P^1$ &
$\Z_2$ \\
\hline
$\Z_2$ $(A_1)$ &
$\left(\begin{array}{rrr} -1 & 0 & 0 \\ 0 & -1 & 0 \\ 0 & 0 & 1 \end{array} \right)$ &
$\begin{array}{c} 2e_1 = 0 \\ 2e_2 = 0 \\ \end{array}$ &
$16$ &
$\Z_2\times \Z_2$ &
$6 \times \P^1$ &
$\Z_2$ \\
\hline
$\Z_3$ $(A_2)$ &
$\left(\begin{array}{rrr} 0 & 0 & 1 \\ 1 & 0 & 0 \\ 0 & 1 & 0 \end{array} \right)$ &
$\begin{array}{c} e_1 = e_3 \\ e_1 = e_2 \\ \end{array}$ &
$1$ &
$\Z_2$ &
$1 \times \P^1$ &
$\Z_2$ \\
\hline
$\Z_4$ $(A_3)$ &
$\left(\begin{array}{rrr} 0 & -1 & 0 \\ 1 & 0 & 0 \\ 0 & 0 & 1 \end{array} \right)$ &
$\begin{array}{c} e_1 = e_2 \\ 2e_1 = 0 \\ \end{array}$ &
$4$ &
$\Z_2$ &
$4 \times \P^1$ &
$\Z_2$ \\
\hline
\end{tabular}}
\vspace{0.25cm}

The set of points with non-cyclic isotropy consists of $24$ points for $D_4$, which become $4$ points in the quotient. Next, there are $36$ points with isotropy $D_6$ in $3$ families associated to subgroups. Image of this set contains $12$ points. There are also $4$ points fixed by $S_4$.

Each curve for the first $\Z_2$ class in the table contains $3$ points with isotropy $D_8$, not identified in the quotient, and $1$ point fixed by $S_4$. Each curve for the second $\Z_2$ class contains $2$ points with isotropy $D_4$ and $2$ with isotropy $D_8$, mapped to different points of $Y$. All fixed points of $S_4$ lie on the $\P^1$ for $\Z_3$. And each curve for $\Z_4$ there are $3$ points with isotropy $D_8$, not identified in $Y$, and $1$ of $S_4$.

The virtual Poincar\'e polynomials of the strata are the following:
\begin{align*}
P_3(t) &= P_Y(t) - 15(1 + t^2 - 4) - 20 =  t^6 + t^4 + 4t^3 - 14t^2 + 26,\\
P_2(t) &= 10(1 + t^2 - 4)(1 + t^2) + \mu_0((1 + \eps \cdot 2t + t^2 - 4)(1 + (1 + \eps)t^2)) +\\ &\quad + 4\mu_0((1 + \eps \cdot 2t + t^2 - 4)(1 + (2 + \eps)t^2)) = 19t^4 + 10t^3 - 42t^2 - 45,\\
P_1(t) &= 16(1 + 4t^2) + 4(1 + 3t^2) = 76t^2 + 20,
\end{align*}
and the Poincar\'e polynomial of $X$ is $$P_X(t) = t^6 + 20t^4 + 14t^3 + 20t^2 + 1.$$

\subsection{Summary}\label{summary}

The following theorem summarizes the results of this paper. Its proof consists of computations presented in sections \ref{examples} and \ref{results}.

\begin{thm}\label{last_theorem}
The complete list of Poincar\'e polynomials of Kummer 3-folds is given in table 2.
\end{thm}

\begin{table}[!hp]
\renewcommand{\arraystretch}{2}
\vspace{0.25cm}
\centerline{\begin{tabular}{ccc}
\textbf{group} & \textbf{section} & \textbf{polynomial} \\
\hline
$D_4(1)$ &
\ref{case_D4_1} &
$t^6 + 51t^4 + 8t^3 + 51t^2 + 1$ \\
$D_4(2)$ &
\ref{case_D4_2} &
$t^6 + 21t^4 + 20t^3 + 21t^2 + 1$ \\
$D_4(3)$ &
\ref{case_D4_3} &
$t^6 + 15t^4 + 8t^3 + 15t^2 + 1$ \\
$D_4(4)$ &
\ref{case_D4_4} &
$t^6 + 15t^4 + 8t^3 + 15t^2 + 1$ \\
\hline
$D_6(1)$ &
\ref{case_D6_1} &
$t^6 + 15t^4 + 32t^3 + 15t^2 + 1$ \\
$D_6(2)$ &
\ref{case_D6_2} &
$t^6 + 15t^4 + 32t^3 + 15t^2 + 1$ \\
$D_6(3)$ &
\ref{case_D6_3} &
$t^6 + 7t^4 + 16t^3 + 7t^2 + 1$ \\
\hline
$D_8(1)$ &
\ref{case_D8_1} &
$t^6 + 36t^4 + 14t^3 + 36t^2 + 1$ \\
$D_8(2)$ &
\ref{case_D8_2} &
$t^6 + 15t^4 + 8t^3 + 15t^2 + 1$ \\
\hline
$D_{12}$ &
\ref{case_D12} &
$t^6 + 21t^4 + 20t^3 + 21t^2 + 1$ \\
\hline
$A_4(1)$ &
\ref{case_A4_1} &
$t^6 + 19t^4 + 8t^3 + 19t^2 + 1$ \\
$A_4(2)$ &
\ref{case_A4_2} &
$t^6 + 7t^4 + 8t^3 + 7t^2 + 1$ \\
$A_4(3)$ &
\ref{case_A4_3} &
$t^6 + 7t^4 + 8t^3 + 7t^2 + 1$ \\
\hline
$S_4(1)$ &
\ref{case_S4_1} &
$t^6 + 20t^4 + 14t^3 + 20t^2 + 1$ \\
$S_4(2)$ &
\ref{case_S4_2} &
$t^6 + 11t^4 + 8t^3 + 11t^2 + 1$ \\
$S_4(3)$ &
\ref{case_S4_3} &
$t^6 + 11t^4 + 8t^3 + 11t^2 + 1$ \\
\end{tabular}}
\vspace{0.25cm}
\renewcommand{\arraystretch}{1}
\caption{Poincar\'e polynomials of Kummer 3-folds}
\end{table}

Poincar\'e polynomials are not sufficient to distinguish varieties obtained by the $3$-dimensional Kummer construction. The question whether the Kummer $3$-folds with equal Poincar\'e are isomorphic has not been solved yet. We only know that subgroups of $SL(3,\Z)$ which are not $\Z$-equivalent define different structures of the quotients of $A^3$ by their actions (in \ref{comment_D4} we discuss two most similar cases). However, this suggests only that even if the isomorphism between Kummer $3$-folds exists, it does not come from the construction in a natural way. Possibly some other invariants work better than Poincar\'e polynomials in this problem. The next step is understanding the cone of effective divisors and the cone of curves of Kummer $3$-folds.

The second point is the relation between Kummer $3$-folds obtained for $\Z$-classes which are dual in the sense of definition \ref{def_duality}. Comparing the list of pairs of dual $\Z$-classes (proposition \ref{dual_pairs}) and the results of computations we can see that Kummer varieties constructed from dual $\Z$-classes have equal Poincar\'e polynomials. It can be investigated whether this is specific for $3$-dimensional Kummer varieties, or is still true in higher dimensions. If it turns out that varieties obtained from the action of dual groups are in fact isomorphic, we can ask whether in higher dimensions the construction for dual groups also give isomorphic varieties.

One more idea for further investigation of Kummer $3$-folds is to look at inclusion of groups and induced rational maps of varieties. We compared the results of computations and the diagram \ref{groups_relations_diagram} of $\Z$-classes inclusions. There are few pairs of groups which are not dual but lead to Kummer varieties with equal Poincar\'e polynomials: $D_4(2)$ and $D_{12}$, $D_4(3)$ and $D_8(2)$, $D_4(4)$ and $D_8(2)$. Note that all these pairs are inclusions up to $\Z$-equivalence. We are interested in finding any significant consequences or generalizations of this observation. It can be also checked whether Kummer $3$-folds appear on the lists of known examples of Calabi--Yau varieties.

\bibliographystyle{amsalpha}
\bibliography{xbib}

\end{document}